\documentclass[11pt]{article}

\usepackage[margin=1in]{geometry}
\usepackage{amsmath,amssymb,amsthm,mathtools,bm}
\usepackage{graphicx,booktabs}
\usepackage{placeins}
\usepackage{algorithm}
\usepackage{algpseudocode}
\usepackage{xcolor}
\usepackage[colorlinks=true,allcolors=blue]{hyperref}
\usepackage[capitalise,noabbrev]{cleveref}

\newcommand{\R}{\mathbb{R}}
\newcommand{\E}{\mathbb{E}}
\newcommand{\dd}{\,\mathrm{d}}
\newcommand{\KL}{\mathrm{D}_{\mathrm{KL}}}

\newcommand{\cc}{c}
\newcommand{\Sc}{S}
\newcommand{\lw}{\lambda}

\newcommand{\refeta}{\eta}
\newcommand{\veta}{v_{\eta}}
\newcommand{\uth}{u_{\theta}}

\newcommand{\Ghat}{\widehat{G}}
\newcommand{\meas}{\nu}

\newcommand{\lab}{z}
\newcommand{\lamres}{\alpha_{\mathrm{res}}}

\newtheorem{theorem}{Theorem}[section]
\newtheorem{proposition}[theorem]{Proposition}
\newtheorem{lemma}[theorem]{Lemma}

\theoremstyle{definition}
\newtheorem{assumption}[theorem]{Assumption}
\newtheorem{definition}[theorem]{Definition}

\title{Score-based neural ordinary differential equations for reaction-diffusions}
\author{Jianguo Hou\thanks{Email: \href{mailto:jianguo@email.sc.edu}{\texttt{jianguo@email.sc.edu}}.}, Wuchen Li\thanks{Corresponding author. Email: \href{mailto:wuchen@mailbox.sc.edu}{\texttt{wuchen@mailbox.sc.edu}}.}\\[4pt]
\small Department of Mathematics\\
\small University of South Carolina\\
\small Columbia, SC 29208, USA}
\date{}

\begin{document}
\maketitle

\begin{abstract}
A diffusion equation can be rewritten as a transport equation driven by the score function of the evolving density. Conventional probability-flow formulations, which move particles deterministically with this score-driven velocity, typically conserve total mass. In this work, we develop a Lagrangian method for reaction-diffusion equations in which each particle carries a position and a physical weight. The characteristic ODEs incorporate source terms into the evolution of the log-density and its gradient, while particle weights track changes in mass. We combine this representation with a learned velocity field and weighted particle training to approximate these dynamics. For non-increasing reaction rates, a relative-entropy inequality supports the mass-normalized self-consistency residual as the training objective. This objective can control the error of the induced density when the growth rate is bounded. The algorithm divides the time interval into short windows. In each window, a frozen velocity generates the particle data for a regression of the network velocity. We draw the particles from period-adapted initial sampling distributions and weight them in the loss partly by their predicted future mass, because the reaction can make regions with little initial mass important later. We numerically investigate the method on the Ornstein--Uhlenbeck (OU) equation with a linear source term, Fisher--KPP equations, and Wasserstein--Fisher--Rao gradient flows. 
\end{abstract}

\section{Introduction}\label{sec:intro}

Reaction-diffusion-advection equations
\[
  \partial_t\rho+\nabla\!\cdot(\rho\,b)=\gamma\Delta\rho+g(t,x,\rho)\,\rho
\]
describe a density $\rho:[0,T]\times\R^d\to(0,\infty)$ that is moved by a drift $b:[0,T]\times\R^d\to\R^d$, spread by diffusion with a constant $\gamma>0$ and created or removed locally at the rate $g:[0,T]\times\R^d\times(0,\infty)\to\R$. Such equations model the spread of populations in ecology~\cite{cantrell2003spatial}, the geographic spread of epidemics~\cite{noble1974plagues}, the invasion of tumor cells~\cite{swanson2003virtual} and many other processes in biology~\cite{murray2002mathbio}. For example, the Fisher--KPP equation~\cite{fisher1937wave,kpp1937} is the basic example in population dynamics with a motility and a proliferation rate. The same structure appears in Wasserstein--Fisher--Rao gradient flows~\cite{chizat2018interpolating} of the Kullback--Leibler (KL) divergence, where the reaction changes the relative mass of distant regions, as in birth-death sampling~\cite{lu2019birthdeath}. In these problems the evolution path of the density matters: when a population grows, and when and where a new mode appears. Because of the reaction, the total mass is in general not conserved, and the density in a region can grow by orders of magnitude over the time interval. Grid-based solvers handle such equations well in low dimensions, but their cost grows exponentially with the dimension, which motivates particle methods that represent the density without a spatial grid.

Particle methods for such equations face several difficulties. First, the reaction can increase or decrease the weights of the particles along a trajectory. Therefore, few particles may gradually come to carry most of the mass, and the effective sample size collapses. Stochastic methods restore the balance by resampling or branching, but introduce additional noise~\cite{delmoral2004feynman,delmoral2006smc}. Second, the density is needed at the particle positions when the rate $g$ depends on $\rho$, as for the Fisher--KPP equation. We have to estimate the density if particles only store positions and weights, for example on histogram bins. Such estimates become noisy and biased as the dimension increases. Furthermore, to move the particles deterministically, we need $\nabla\log\rho$ to rewrite the diffusion term, which typically adds a further estimation step.

Deterministic probability-flow methods move the particles with the velocity $b-\gamma\nabla\log\rho$, where $\nabla\log\rho$ is the score~\cite{song2021sde,maoutsa2020interacting,boffi2023probability}. As solvers of a given equation, they have mainly been developed for mass-conserving equations. In this paper we extend the deterministic probability-flow approach to equations with the local reaction term. Given a velocity field, each particle carries its position, log-density, score and log-weight. The reaction enters the characteristic system directly (Proposition~\ref{prop:carry}). The velocity is represented by a neural network and trained to be self-consistent. With this parameterization the characteristic system is a score-based neural ODE. Unlike inference from snapshot data, we solve a forward problem: no samples of $\rho_t$ are available, and the velocity is determined by self-consistency. Compared with stochastic particle methods, our particles move deterministically and carry their own density, while the velocity must be learned and the trajectories describe the density rather than the sample paths of the stochastic process.

To use finite particles, we need to answer the question of where to place the particles and how to weight them in the training. In the mass-conserving case, adaptive sampling addresses the placement~\cite{Liao2025apfrm}. With a reaction term, a region that holds a large part of the mass in the future may start from very low initial density. This region may not have enough particles and almost no weight in training if we start from the physical initial distribution. Therefore, we draw the particles from different initial sampling distributions. Multiple importance sampling and the design of proposal distributions for particle filters and sequential Monte Carlo have been studied in~\cite{owen2000safe,cornuet2012amis,whiteley2014twisted,heng2020controlled}. We also adaptively balance the weights of particles in the loss, which gives regions that become important later more influence on the fit.

We train the velocity window by window. In each window, a frozen velocity generates the particle data, and the network is fitted to these fixed data. The training loss combines a regression of the network velocity onto the velocity given by the carried score with a regularizer from the PDE satisfied by the exact velocity. The main contributions of our work are:
\begin{itemize}
  \item a closed characteristic system for reaction-diffusion-advection equations, with the reaction evaluated at the carried density (Section~\ref{sec:weighted_char});
  \item a relative-entropy inequality, which motivates the mass-normalized self-consistency residual as the objective and bounds the error linearly when the growth rate is bounded (Section~\ref{sec:kl_control});
  \item period-adapted initial sampling distributions and adaptive regression weights for training (Section~\ref{sec:banks});
  \item a training algorithm in which each iteration is a regression on labels produced by a frozen velocity (Section~\ref{sec:scheme});
  \item tests on Ornstein--Uhlenbeck (OU) equations with a linear source, Fisher--KPP equations and Wasserstein--Fisher--Rao gradient flow problems, evaluated against reference solutions along the path (Section~\ref{sec:numerics}).
\end{itemize}

\paragraph{Related work.} Several other particle methods address these difficulties for related equations. Stochastic particle methods keep the random sample paths. A linear reaction term can be represented by weights along these paths through the Feynman--Kac formula and by interacting particle systems with selection or branching~\cite{delmoral2004feynman}. When the rate depends on the density, the density is estimated from the particles~\cite{hu2026sgip}. Two recent solvers treat reaction-diffusion-advection equations with randomly moving particles: the stochastic genetic interacting particle method~\cite{hu2026sgip} integrates the reaction on histogram bins with adaptive resampling, and Lyu and Lei~\cite{lyu2025branching} use a branching birth-death step. In one dimension, Mascagni~\cite{mascagni1995deterministic} derives deterministic ODEs for the positions of particles.

 Deterministic particle methods for diffusion include~\cite{degond1990deterministic,carrillo2019blob}; Maoutsa et al.~\cite{maoutsa2020interacting} estimate the gradient of the log-density from the particles, and Willems et al.~\cite{willems2026sformulation} let each particle carry $\log\rho$ and compute its derivatives from neighboring particles. Boffi and Vanden-Eijnden~\cite{boffi2023probability} show that the score obeys an ODE along the trajectories and that a global objective comparing the learned score with this carried score can bound the KL divergence, while their algorithm estimates the score at each time step from the current particles by score matching; see also~\cite{yang2022inexactlangevin,chen2023probflowfast} for sampling with an estimated score. Shen et al.~\cite{shen2022selfconsistency} formulated the self-consistency of the Fokker--Planck equation, and Shen and Wang~\cite{shenwang2023einn} and Lu et al.~\cite{lu2024mfsbtm} control the KL divergence for mean-field equations by density-weighted residuals of the velocity and of the score. With Neural ODEs, the velocity is parameterized by a neural network and integrated along the trajectories~\cite{chen2018neuralode}. Li et al.~\cite{li2023scvm} fit a neural velocity iteratively. Zhou et al.~\cite{zhou2024scoreode,zhou2025mfc} parameterize the velocity, propagate the density and the score along its trajectories by ODEs, and minimize a mean-field control objective. Li and Liu~\cite{li2026neuraltransport} carry the density and a generalized score for aggregation-diffusion equations. Liu et al.~\cite{Zhang2023pinf} fit a log-density network to values integrated along the characteristics. Wu and Liao~\cite{Liao2025apfrm} obtain the log-density from a continuous normalizing flow and minimize the probability flow residual on adaptively sampled points. They bound the KL divergence by this residual weighted by the estimated density. Other neural solvers include physics-informed neural networks and the deep Galerkin method, which minimize the PDE residual at collocation points~\cite{raissi2019pinn,sirignano2018dgm}, and Deep BSDE~\cite{han2018solving} and deep splitting~\cite{beck2021deepsplitting} methods for nonlinear parabolic equations in high dimensions. Score-PINN learns the score of Fokker--Planck solutions~\cite{hu2024scorepinn}.

On finite positive measures, transport combined with the change of mass is described by the Wasserstein--Fisher--Rao (WFR) or Hellinger--Kantorovich geometry~\cite{chizat2018interpolating,kondratyev2016newdistance,liero2018optimal}, and the equation of Section~\ref{sec:wfr} is a WFR gradient flow. Reaction-diffusion equations also have gradient-flow formulations with generalized optimal transport metrics~\cite{li2022reactiondiffusion}, generalized JKO schemes in these metrics~\cite{fu2023highorder}, and splitting schemes~\cite{gallouet2017jko,gallouet2019ubotsplitting} based on the JKO scheme~\cite{jko1998}. Seo et al.~\cite{seo2026hkjko} develop a neural JKO scheme for WFR gradient flows, whose form includes~\eqref{eq:klwfr_pde}. Birth-death dynamics, sequential Monte Carlo and reweighted normalizing flows use WFR-type geometries to sample from a known target~\cite{lu2019birthdeath,crucinio2025smcwfr,duan2026wfrjko}. Some other work fits a velocity and a growth rate $g(x,t)$ based on snapshot data~\cite{sha2024tigon,zhang2025upfi,sha2025ruot,peng2026wfrfm}.

Section~\ref{sec:method} derives the characteristic system, the residual inequality and the continuous objective. Section~\ref{sec:scheme} gives the algorithm. Section~\ref{sec:numerics} reports the experiments, and Section~\ref{sec:conclusion} summarizes the results and discusses the limitations and future work. Appendix~\ref{app:proofs} derives the characteristic system and the weighted identities, and Appendix~\ref{app:entropy} proves the residual inequality. Appendix~\ref{app:impl} describes the implementation, and Appendix~\ref{app:evaluation} the reference solutions and the evaluation.

\section{Self-consistent characteristic flows}\label{sec:method}

We start from the reaction-diffusion-advection equation, and derive the ODEs carried by each particle for a given velocity (Sections~\ref{sec:carry} and~\ref{sec:weighted_char}). We then show that the self-consistency residual controls the error of the induced density (Section~\ref{sec:kl_control}), evaluate this residual with the carried score (Section~\ref{sec:carried_score}), and build the continuous objective (Section~\ref{sec:objective}).

\subsection{Problem setting}\label{sec:carry}

We consider solving the following reaction-diffusion-advection equations:
\begin{equation}\label{eq:rd_general}
  \partial_t \rho
  \;+\;
  \nabla\!\cdot\!(\rho\, b)
  \;=\;
  \gamma\,\Delta\rho
  \;+\;
  g(t,x,\rho)\,\rho,
  \qquad
  \text{given} \quad\rho_0,
\end{equation}
where $b=b(t,x)$ is the drift term, $\gamma\,\Delta\rho$ is the diffusion term with a constant $\gamma>0$ and $g(t,x,\rho)\,\rho$ is the multiplicative reaction term. We assume that $\rho_0>0$ has finite mass
$M_0=\int\rho_0(x)\,\mathrm{d}x$. For a smooth, strictly positive density $\rho$, we can rewrite \eqref{eq:rd_general} as a continuity equation with a source term using $\Delta\rho = \nabla\!\cdot\!(\rho\,\nabla\log\rho)$:
\begin{equation}\label{eq:rd_transport}
  \partial_t\rho
  \;+\;
  \nabla\!\cdot\!(\rho\,v_\rho)
  \;=\;
  g(t,x,\rho)\,\rho,
  \qquad
  v_\rho \;:=\; b - \gamma\,\nabla\log\rho.
\end{equation}
Different from the transport equation, the multiplicative reaction term $g$ will increase or decrease the density locally. As a result, the conservation of total mass often does not hold.
We further make the following assumptions on the problem:

\begin{assumption}\label{ass:continuum}
The drift $b$ is smooth on $[0,T]\times\R^d$.
The initial density $\rho_0$ is smooth, strictly positive and of finite mass, and \eqref{eq:rd_general} has a smooth, strictly positive solution $\rho$ on $[0,T]$, and no other strictly positive solution that is $C^2$ in $x$ and piecewise $C^1$ in $t$. The reaction rate $g(t,x,r)$ is continuous, and twice continuously differentiable in $(x,r)$ for $r>0$, and for every $(t,x)$ the map $\ell\mapsto g(t,x,e^{\ell})$ is locally Lipschitz and non-increasing.
\end{assumption}

\subsection{Weighted characteristics}\label{sec:weighted_char}

We next represent and approximate the solution of \eqref{eq:rd_general} by weighted particles moving along the characteristics of \eqref{eq:rd_transport}. To get access to $v_\rho$, we first define the admissibility conditions for a candidate velocity $v$. Then we will show that an admissible $v$ is exactly $v_\rho$ if it satisfies the self-consistency condition in Proposition~\ref{prop:selfconsistent}.
\begin{definition}[Admissible velocity]\label{def:admissible}
A velocity field $v:[0,T]\times\R^d\to\R^d$ is \emph{admissible} if there are times $0=T_0<T_1<\dots<T_{L}=T$ such that on each window the restriction of $v$ to $(T_{l},T_{l+1})\times\R^d$ and its spatial derivatives up to order three extend continuously in $(t,x)$ to $[T_{l},T_{l+1}]\times\R^d$, and $\nabla v$ is bounded.
\end{definition}
Now consider the following transport-reaction equation, where the velocity $v$ is given and admissible:
\begin{equation}\label{eq:mu_pde}
  \partial_t\mu+\nabla\!\cdot(\mu\,v)=g(t,x,\mu)\,\mu,\qquad \mu(0)=\rho_0 .
\end{equation}
For this initial value problem, we can show that the solution exists and is determined by an ODE system:

\begin{proposition}\label{prop:carry}
Let Assumption~\ref{ass:continuum} hold, let $v$ be admissible (Definition~\ref{def:admissible}), and write $\nabla v$ for its Jacobian. From initial distribution $X_0\sim\rho_0/M_0$, $v$ drives the following ODE system, which we call the characteristic system
\begin{align}
  \dot X_t &= v(t,X_t),\label{eq:char_x}\\
  \dot\ell_t &= -\nabla\!\cdot v(t,X_t)+g\bigl(t,X_t,e^{\ell_t}\bigr),\label{eq:char_ell}\\
  \dot\Sc_t &= -\bigl(\nabla v(t,X_t)\bigr)^{\top}\Sc_t-\nabla\bigl(\nabla\!\cdot v\bigr)(t,X_t)+\nabla_x g\bigl(t,X_t,e^{\ell_t}\bigr),\label{eq:char_score}\\
  \dot\lw_t &= g\bigl(t,X_t,e^{\ell_t}\bigr),\label{eq:char_lw}
\end{align}
with initial values $\ell_0=\log\rho_0(X_0)$, $\Sc_0=\nabla\log\rho_0(X_0)$ and $\lw_0=\log M_0$. 
Here $\nabla_x g$ is written in terms of the particle state,
\begin{equation}\label{eq:nabla_g}
  \nabla_x g\bigl(t,X_t,e^{\ell_t}\bigr)=(\partial_x g)\bigl(t,X_t,e^{\ell_t}\bigr)+\partial_\rho g\bigl(t,X_t,e^{\ell_t}\bigr)\,e^{\ell_t}\,\Sc_t ,
\end{equation}
where $\partial_x g$ and $\partial_\rho g$ are the partial derivatives of $g(t,x,r)$ with respect to $x$ and $r$. Then, through this system, $v$ induces a density $\mu=\mu[v]$, which is the solution of \eqref{eq:mu_pde}. Moreover,
\begin{enumerate}
  \item[(i)] along every trajectory of $X_t$, $\ell_t=\log\mu(t,X_t)$ and $\Sc_t=\nabla\log\mu(t,X_t)$, and \eqref{eq:nabla_g} is the total spatial gradient of $x\mapsto g(t,x,\mu(t,x))$;
  \item[(ii)] $\mu(t,\cdot)$ is the pushforward density of $\rho_0/M_0$ under the flow map $X_0\mapsto X_t$, with each trajectory weighted by $e^{\lw_t}$,
  \begin{equation}\label{eq:mu_def}
    \int_A\mu(t,x)\,\mathrm{d}x=\mathbb{E}_{X_0\sim\rho_0/M_0}\bigl[e^{\lw_t}\mathbf 1\{X_t\in A\}\bigr]\qquad\text{for every Borel set }A.
  \end{equation}
\end{enumerate}
\end{proposition}

Furthermore, the induced density $\mu[v]$ can match the solution of $\eqref{eq:rd_general}$ under the following condition:

\begin{proposition}[Self-consistent velocity]\label{prop:selfconsistent}
Let $\mu=\mu[v]$ be the density induced by $v$ in Proposition~\ref{prop:carry}. Then $\mu$ is the solution of \eqref{eq:rd_general} if
\begin{equation}\label{eq:self_consistency}
  v=b-\gamma\nabla\log\mu\qquad\text{on }[0,T].
\end{equation}
In particular, $v=v_\rho$, $\ell_t=\log\rho(t,X_t)$ and $\Sc_t=\nabla\log\rho(t,X_t)$ along every trajectory.
Conversely, if the exact velocity $v_\rho=b-\gamma\nabla\log\rho$ is admissible, it induces $\mu[v_\rho]=\rho$ and therefore satisfies \eqref{eq:self_consistency}.
\end{proposition}

We therefore let the velocity $v$ approach $b-\gamma\nabla\log\mu[v]$ to match condition $\eqref{eq:self_consistency}$.

\subsection{Stability of the induced density}\label{sec:kl_control}

We now measure how far a velocity is from self-consistency. Let $v$ be admissible and let $\mu:=\mu[v]$ be its induced density. We define the self-consistency residual and its mass-normalized mean square, integrated over time,
\begin{equation}\label{eq:residual}
  R_v := v-\bigl(b-\gamma\nabla\log\mu\bigr),
  \qquad
  J_{\mathrm{reg}}(v):=\int_0^T\frac{1}{M_\mu(t)}\int\mu\,|R_{v}|^2\,\mathrm{d}x\,\mathrm{d}t,
\end{equation}
where $M_\mu(t):=\int\mu(t,x)\dd x$, which we assume to be finite on $[0,T]$. By Proposition~\ref{prop:selfconsistent}, $J_{\mathrm{reg}}(v)=0$ means that $\mu$ is the solution of \eqref{eq:rd_general}. We normalize by the mass so that every time contributes on the same scale when the mass changes. Writing $\mu_t(\cdot)=\mu(t,\cdot)$ and $\rho_t(\cdot)=\rho(t,\cdot)$, we measure the distance between $\mu_t$ and $\rho_t$ by the generalized KL divergence for finite positive densities:
\begin{equation}\label{eq:gen_kl}
  \KL(\mu_t\|\rho_t)
  :=\int\rho_t\,\phi\Bigl(\frac{\mu_t}{\rho_t}\Bigr)\mathrm dx,
  \qquad
  \phi(s):=s\log s-s+1\;\ge0 ,
\end{equation}
which vanishes only if $\mu_t=\rho_t$.
\begin{theorem}[Residual inequality]\label{thm:stability}
Let Assumption~\ref{ass:continuum} hold, let $v$ be admissible, with induced density $\mu=\mu[v]$, and assume that $\rho$ and $\mu$ decay fast enough at infinity (Assumption~\ref{ass:decay}).
\begin{enumerate}
  \item[(i)] For almost every $t\in[0,T]$,
  \begin{equation}\label{eq:kl_rate}
    \frac{\mathrm d}{\mathrm dt}\KL(\mu_t\|\rho_t)
    \le\frac1{4\gamma}\int\mu_t\,|R_v|^2\,\mathrm dx
       +\int g(t,x,\rho_t)\,\rho_t\,\phi(\mu_t/\rho_t)\,\mathrm dx .
  \end{equation}
  \item[(ii)] Let $G(t):=\sup_x[g(t,x,\rho(t,x))]_+$ be the largest growth rate of the solution at time $t$. If $G\in L^1(0,T)$, then
  \begin{equation}\label{eq:kl_bound}
    \sup_{0\le t\le T}\KL(\mu_t\|\rho_t)\le\frac{C_G}{4\gamma}\,J_{\mathrm{reg}}(v),
    \qquad
    C_G:=\sup_{0\le s\le t\le T}M_\mu(s)\,e^{\int_s^tG(\tau)\,\mathrm d\tau}.
  \end{equation}
\end{enumerate}
\end{theorem}

The argument follows the relative-entropy estimates for mass-conserving equations~\cite{boffi2023probability,shenwang2023einn,lu2024mfsbtm,Liao2025apfrm}. With a reaction, the growth rate adds the second term of \eqref{eq:kl_rate}, and the monotonicity of $g$ removes the comparison of the two reaction rates (Appendix~\ref{app:inequality}).

This motivates us to minimize $J_{\mathrm{reg}}$. Evaluating $R_v$ requires the score $\nabla\log\mu$ of the density induced by $v$. We will use the score carried along the trajectory to compute the $J_{\mathrm{reg}}$.
\subsection{Carried score and the weighted residual}\label{sec:carried_score}
System~\eqref{eq:char_x}-\eqref{eq:char_score} is closed: given the particle state $(X_t, \ell_t, \Sc_t)$ and the velocity field $v$, all derivatives can be evaluated. We can therefore propagate an ensemble of particles forward in time, with each particle carrying the log-density and the score of the induced density $\mu[v]$ instead of estimating them from the particles. Moreover, the initial distribution of the particles can differ from $\rho_0/M_0$.

\begin{lemma}[Weighted residual identity]\label{lem:weighted_residual}
Let Assumption~\ref{ass:continuum} hold, $\mu=\mu[v]$ be the solution of \eqref{eq:mu_pde}, and
$q_0$ be a probability density with $q_0>0$ wherever $\rho_0>0$. Draw $X_0\sim q_0$, set
\[
  \ell_0=\log\rho_0(X_0),\qquad \Sc_0=\nabla\log\rho_0(X_0),\qquad \lw_0=\log\rho_0(X_0)-\log q_0(X_0),
\]
and propagate $(X_t,\ell_t,\Sc_t,\lw_t)$ by \eqref{eq:char_x}--\eqref{eq:char_lw}. Then, with $\mathbb{E}$ taken over $X_0\sim q_0$:
\begin{enumerate}
  \item[(i)] along every trajectory, $\ell_t=\log\mu(t,X_t)$ and $\Sc_t=\nabla\log\mu(t,X_t)$;
  \item[(ii)] for every $t\in[0,T]$ and every measurable $f$ that is non-negative or satisfies $\int\mu|f|\,\mathrm{d}x<\infty$,
  \begin{equation}\label{eq:weighted_identity}
    \mathbb{E}\bigl[e^{\lw_t}f(t,X_t)\bigr]=\int\mu(t,x)\,f(t,x)\,\mathrm{d}x ;
  \end{equation}
  \item[(iii)] in particular, we have
  \begin{equation}\label{eq:residual_identity}
  \begin{aligned}
    M_\mu(t)&=\int\mu(t,x)\,\mathrm{d}x=\mathbb{E}[e^{\lw_t}],\\
    \int\mu\,|R_v|^2\,\mathrm{d}x
    &=\mathbb{E}\Bigl[e^{\lw_t}\,\bigl|v(t,X_t)-b(t,X_t)+\gamma\Sc_t\bigr|^2\Bigr].
  \end{aligned}
  \end{equation}
\end{enumerate}
Proposition~\ref{prop:carry} is the special case $q_0=\rho_0/M_0$, for which $\lw_0=\log M_0$.
\end{lemma}

According to this lemma, we can choose any probability density $q_0$ that is positive wherever $\rho_0>0$ as the initial sampling distribution, as long as we change the initial log-weight $\lw_0$ correspondingly; in the terminology of importance sampling, $q_0$ is the proposal distribution and $e^{\lw_0}=\rho_0/q_0$ is the importance weight. Starting from initial distribution $q_0$, for any admissible velocity field $v$ (Definition~\ref{def:admissible}) to drive the characteristic system, we can evaluate the objective functional $J_{\mathrm{reg}}$ of \eqref{eq:residual} using the carried score $\Sc_t$:

\begin{equation}\label{eq:J_reg}
  J_{\mathrm{reg}}(v)
  =\int_0^T\frac{\mathbb{E}\bigl[e^{\lw_t}\,|v(t,X_t)-b(t,X_t)+\gamma \Sc_t|^2\bigr]}{\mathbb{E}[e^{\lw_t}]}\,\mathrm{d}t ,
\end{equation}
where $\mathbb{E}$ is taken over $X_0\sim q_0$ and the initial values are those of Lemma~\ref{lem:weighted_residual}.

\subsection{The continuous objective}\label{sec:objective}
The second order derivatives $\nabla(\nabla\!\cdot v)$ do not appear explicitly in \eqref{eq:J_reg}. They affect $J_{\mathrm{reg}}$ only through the carried score $\Sc_t$, which integrates them along the trajectory by \eqref{eq:char_score}. To regularize the training, we consider the equation satisfied by the exact velocity, which explicitly contains $\nabla(\nabla\!\cdot v)$. Starting from \eqref{eq:rd_transport}, $v_\rho=b-\gamma\nabla\log\rho$ satisfies

\begin{equation}\label{eq:velocity_equation}
  \partial_t v_\rho
  =\partial_t b+\gamma\nabla(\nabla\!\cdot v_\rho)+\nabla\bigl(v_\rho\cdot(b-v_\rho)\bigr)
   -\gamma\,\partial_x g(t,x,\rho)-\rho\,\partial_\rho g(t,x,\rho)\,(b-v_\rho).
\end{equation}
Based on \eqref{eq:velocity_equation}, for an admissible $v$ that is continuously differentiable in $t$ on each window, we can define the residual $P(v)$ as follows:
\begin{equation}\label{eq:velocity_pde_residual}
  P(v):=\partial_t v-\partial_t b-\gamma\nabla(\nabla\!\cdot v)-\nabla\bigl(v\cdot(b-v)\bigr)+\gamma\,\partial_x g(t,x,\mu)+\mu\,\partial_\rho g(t,x,\mu)\,(b-v),
\end{equation}
where $\mu=\mu[v]$ is the density induced by $v$. The residual $P(v)$ is obtained from \eqref{eq:velocity_equation} by replacing $v_\rho$ with $v$ and $\rho$ with $\mu[v]$, so it depends only on $v$. If $v_\rho$ is admissible, then $\mu[v_\rho]=\rho$ by Proposition~\ref{prop:selfconsistent}, and $P(v_\rho)=0$ is exactly \eqref{eq:velocity_equation}. The derivations of \eqref{eq:velocity_equation} and of the following identity are given in Appendix~\ref{app:proof_velocity}.

\begin{proposition}[Velocity equation and the self-consistency residual]\label{prop:velocity}
Under Assumption~\ref{ass:continuum}, for every admissible $v$ that is continuously differentiable in $t$ on each window,
\begin{equation}\label{eq:residual_transport}
  \partial_t R_v+\nabla\bigl(v\cdot R_v\bigr)-\mu\,\partial_\rho g(t,x,\mu)\,R_v \;=\; P(v).
\end{equation}
Consequently:
\begin{enumerate}
  \item[(i)] if $R_v\equiv0$ on $[0,T]$, then $P(v)\equiv0$ on $[0,T]$;
  \item[(ii)] the exact velocity $v_\rho=b-\gamma\nabla\log\rho$ satisfies the closed equation \eqref{eq:velocity_equation}.
\end{enumerate}
\end{proposition}
We define the regularization objective functional using residual $P(v)$:
\begin{equation}\label{eq:J_res}
  J_{\mathrm{res}}(v):=\int_0^T\frac{1}{M_\mu(t)}\int\mu\,|P(v)|^2\,\mathrm{d}x\,\mathrm{d}t
\end{equation}

Combining $J_{\mathrm{reg}}$ and $J_{\mathrm{res}}(v)$, we get the composite objective functional:
\begin{equation}\label{eq:J_com}
  J_{\mathrm{com}}(v):=J_{\mathrm{reg}}(v)+\lamres\,J_{\mathrm{res}}(v),\quad \lamres\ge0 .
\end{equation}

The term $J_{\mathrm{res}}$ acts only as a regularizer. By \eqref{eq:residual_transport}, $R_v\equiv0$ implies $P(v)\equiv0$, so $J_{\mathrm{com}}$ and $J_{\mathrm{reg}}$ vanish for the same velocities, and since $J_{\mathrm{reg}}\le J_{\mathrm{com}}$, Theorem~\ref{thm:stability}(ii) holds with $J_{\mathrm{com}}$ in place of $J_{\mathrm{reg}}$.

In the next section, we will show how to design a numerical scheme based on the idea of minimizing $J_{\mathrm{com}}(v)$. The numerical scheme is expected to approximate the solution of \eqref{eq:rd_general} by optimizing a training loss motivated by $J_{\mathrm{com}}(v)$. 

\section{Numerical Scheme}\label{sec:scheme}

For any velocity field, a numerical ODE scheme gives a numerical solution of the characteristic system. We represent the velocity field by a neural network. Then we optimize it so that the numerical solution matches the characteristics of the exact velocity.

In the following subsections, we describe how to obtain the numerical ODE solution from a neural network based velocity field $v_\theta$. Then we can capture the score value along the characteristics. We next describe the outer iteration, in which a frozen velocity generates the particle data. Using the information carried by the particles, we construct a training loss motivated by $J_{\mathrm{com}}(v)$ and train the neural network for $v_\theta$ window by window.

\subsection{Time integration and particle data}\label{sec:time}

For a given velocity field $v$, we integrate the characteristic system \eqref{eq:char_x}--\eqref{eq:char_lw} in Proposition~\ref{prop:carry}. To discretize this system, we denote the particle state by $\zeta_t := (X_t,\,\ell_t,\,\Sc_t,\,\lw_t)$ with $X_0\sim q_0$, and write the right-hand side of the characteristic system as $F_\theta(t,\zeta_t)$:
\begin{equation}\label{eq:F_def}
  F_\theta(t,\zeta_t)
  \;=\;
  \begin{bmatrix}
    v_\theta(t,X_t) \\[2pt]
    g(t,X_t,e^{\ell_t}) - \nabla\!\cdot\! v_\theta(t,X_t) \\[2pt]
    -(\nabla v_\theta(t,X_t))^\top \Sc_t - \nabla(\nabla\!\cdot\! v_\theta(t,X_t)) + \nabla_x g(t,X_t,e^{\ell_t}) \\[2pt]
    g(t,X_t,e^{\ell_t})
  \end{bmatrix},
\end{equation}
where the spatial derivatives $\nabla\!\cdot\! v_\theta$, $\nabla(\nabla\!\cdot\! v_\theta)$, and $(\nabla v_\theta)^\top \Sc_t$ are computed by automatic differentiation of the network. We partition $[0,T]$ into the time windows $I_j:=[T_j,T_{j+1}]$, $j=0,\dots,N_w-1$, with $0=T_0<T_1<\dots<T_{N_w}=T$. Each window is divided into $N_j$ uniform subintervals of length $\Delta t = (T_{j+1}-T_j)/N_j$, with $t_{j,n} = T_j+n\Delta t$, $n=0,\dots,N_j$; within a fixed window we write $t_n$. We integrate $\zeta_t$ forward in time by the explicit trapezoidal (Heun) rule; a step from $t_n$ to $t_{n+1}$ has two stages:
\begin{subequations}\label{eq:heun}
\begin{align}
  \text{(Predictor)}\qquad
  \tilde\zeta_{n+1}
  &\;=\;
  \zeta_n + \Delta t\,F_\theta(t_n,\zeta_n),
  \label{eq:heun_pred}\\[2pt]
  \text{(Corrector)}\qquad
  \zeta_{n+1}
  &\;=\;
  \zeta_n
  \;+\;
  \tfrac{\Delta t}{2}\Big[\,F_\theta(t_n,\zeta_n)
  \;+\;
  F_\theta(t_{n+1},\tilde\zeta_{n+1})\,\Big].
  \label{eq:heun_corr}
\end{align}
\end{subequations}

We write $\mathcal L_\theta(t,\zeta):=g(t,X,e^{\ell})-\nabla\!\cdot\! v_\theta(t,X)$ and $\mathcal S_\theta(t,\zeta):=-(\nabla v_\theta(t,X))^\top \Sc - \nabla(\nabla\!\cdot\! v_\theta(t,X)) + \nabla_x g(t,X,e^{\ell})$ for the second and third components of $F_\theta$; the fourth is $g$. For the $i$-th particle, the corrector \eqref{eq:heun_corr} reads componentwise
\begin{align}
  X^i_{n+1}
  &\;=\;
  X^i_n \;+\; \tfrac{\Delta t}{2}\Big[v_\theta(t_n,X^i_n) + v_\theta(t_{n+1},\tilde X^i_{n+1})\Big],
  \label{eq:disc_x}\\[2pt]
  \ell^i_{n+1}
  &\;=\;
  \ell^i_n \;+\; \tfrac{\Delta t}{2}\Big[\mathcal L_\theta(t_n,\zeta^i_n) \;+\; \mathcal L_\theta(t_{n+1},\tilde\zeta^i_{n+1})\Big],
  \label{eq:disc_l}\\[2pt]
  \Sc^i_{n+1}
  &\;=\;
  \Sc^i_n \;+\; \tfrac{\Delta t}{2}\Big[\mathcal S_\theta(t_n,\zeta^i_n) \;+\; \mathcal S_\theta(t_{n+1},\tilde\zeta^i_{n+1})\Big],
  \label{eq:disc_s}\\[2pt]
  \lw^i_{n+1}
  &\;=\;
  \lw^i_n \;+\; \tfrac{\Delta t}{2}\Big[g(t_n,X^i_n,e^{\ell^i_n}) \;+\; g(t_{n+1},\tilde X^i_{n+1},e^{\tilde\ell^i_{n+1}})\Big],
  \label{eq:disc_w}
\end{align}
where $\tilde\zeta^i_{n+1}=(\tilde X^i_{n+1},\tilde\ell^i_{n+1},\tilde\Sc^i_{n+1},\tilde\lw^i_{n+1})$ is the predictor \eqref{eq:heun_pred}. Starting from $N$ particles with $X^i_0\sim q_0$ and the initial values of Lemma~\ref{lem:weighted_residual}, we obtain the particle data $\mathcal D(v_\theta):=\{\zeta^i_n\}_{i,n}$. To make these particles track the solution $\rho$ of \eqref{eq:rd_transport}, we optimize the velocity field $v_\theta$ by minimizing a loss function $L$.

\subsection{Frozen-velocity iteration}\label{sec:iter}
Inside each time window $I_j$, $j=0,\dots,N_w-1$ (see Section~\ref{sec:time}), starting from the initial velocity field $v_{\bar{\theta}^{(0)}}$, we propagate the characteristic system and collect $\mathcal D^{(0)}:=\mathcal D(v_{\bar{\theta}^{(0)}})$. Based on $\mathcal D^{(0)}$, we decrease the loss $L(\theta;\bar{\theta}^{(0)})$ of Section~\ref{sec:frozen} by a fixed number of optimizer steps, which approximately solves the regression problem, and obtain the updated velocity field
\[
  v_{\bar{\theta}^{(1)}} \approx \operatorname*{arg\,min}_{\theta} L(\theta;\bar{\theta}^{(0)}).
\]
We repeat the process until $v_{\bar{\theta}^{(k+1)}}$ is close enough to $v_{\bar{\theta}^{(k)}}$ and the loss $L(\bar{\theta}^{(k+1)};\bar{\theta}^{(k)})$ is small enough.

The labels used in the training depend on the velocity through the characteristic system. Freezing the velocity that generates them makes each outer iteration a weighted regression on fixed particle data, without differentiating through the trajectories. A velocity that reproduces the labels of its own particle data satisfies the self-consistency condition \eqref{eq:self_consistency} at these particles, up to the time-integration error. Since the labels come from the previous velocity, several outer iterations are needed in each window.

\subsection{Loss function}\label{sec:frozen}

We now define the training loss on the particle data generated by a frozen velocity. Within a time window $I_j$, in the $k$-th outer iteration we generate $\mathcal D^{(k)}:=\mathcal D(v_{\bar\theta^{(k)}})$ containing $N$ particles according to the fixed parameter $\bar\theta^{(k)}$. It gives the states $(X^i_n,\ell^i_n,\Sc^i_n,\lw^i_n)$ of particle $i$ at the time points $t_n\in I_j$, $n=0,\dots,N_j-1$. A general formula of the loss function $L$ can be defined using $\mathcal D^{(k)}$:
\begin{equation}\label{eq:J_modular}
  L(\theta;\bar\theta^{(k)})
  :=\frac{1}{N_j}\sum_{n=0}^{N_j-1}\sum_{i=1}^{N}\omega^i_n
   \Bigl[\,\psi_{\rm reg}\bigl(v_\theta(t_n,X^i_n)-v_{\rm label}(t_n,X^i_n)\bigr)
   +\lamres\,\psi_{\rm res}\bigl(P(t_n,X^i_n;v_\theta,\ell^i_n)\bigr)\Bigr],
\end{equation}
where $v_{\rm label}(t_n,X^i_n):=b(t_n,X^i_n)-\gamma\Sc^i_n$ and
\begin{equation}\label{eq:P_theta}
  P(t,x;v_\theta,\ell)
  :=\partial_t v_\theta-\partial_t b-\gamma\nabla(\nabla\!\cdot v_\theta)
    -\nabla\bigl(v_\theta\cdot(b-v_\theta)\bigr)
    +\gamma\,\partial_x g(t,x,e^{\ell})+e^{\ell}\,\partial_\rho g(t,x,e^{\ell})\,(b-v_\theta).
\end{equation}
All derivatives of $v_\theta$ are computed using automatic differentiation. The basic choice is the physical weights with quadratic losses,
\begin{equation}\label{eq:basic_choice}
  \omega^i_n=\frac{e^{\lw^i_n}}{\sum_{i'=1}^Ne^{\lw^{i'}_n}},\qquad \psi_{\rm reg}(\cdot)=\psi_{\rm res}(\cdot)=|\cdot|^2 .
\end{equation}
In all experiments, we replace the squared loss in both terms by the Huber loss~\cite{huber1964robust}, also known as the smooth $L^1$ loss, so that particles with very large residuals do not dominate the loss (Appendix~\ref{app:iteration}). To further improve the accuracy, we also modify the initial distribution $q_0$ and the weights $\omega$, as described next.

\subsubsection{Period-adapted sampling and adaptive weights}\label{sec:banks}

We can get the same $\mu$ starting from different $q_0$ as long as we change $\lw_0=\log\rho_0-\log q_0$ according to Lemma~\ref{lem:weighted_residual}, but we will get different numerical solutions. Notice that we have fewer particles in the region with low density $q_0$ when we use finitely many particles. In addition, the reaction term $g$ changes the weight of a particle along its trajectory by \eqref{eq:char_lw}, starting from $\lw_0=\log\rho_0-\log q_0$.

Therefore, we may not have enough particles around the high density regions in the future if those regions come from the initial low density regions. We need to modify $q_0$ to improve coverage because adding enough particles may lead to high computational cost, especially for high dimensions. To cover the high density regions at different times, we partition $[0,T]$ into three periods $\mathcal{T}_0, \mathcal{T}_1, \mathcal{T}_2$: $\mathcal{T}_0$ for early period, $\mathcal{T}_1$ for middle period and $\mathcal{T}_2$ for late period. For every period $\mathcal T_m$ we construct an initial sampling distribution $q_0^m$ and draw from it a separate ensemble of $N$ particles. We will need three particle ensembles in the training process on the whole interval $[0,T]$. We call the distributions $q_0^m$ the period-adapted sampling distributions. To construct $q_0^m$, we draw pilot particles from a broad initial distribution and integrate their positions and log-weights along the flow of a predicted density $\rho_{\mathrm{pred}}$, a cheap approximation of the solution built from the problem data only (Appendix~\ref{app:banks}). We then build $q_0^m$ mainly from a Gaussian mixture fitted to the initial positions of the pilot particles, each weighted by its share of the mass in period $\mathcal T_m$.

By Lemma~\ref{lem:weighted_residual}, the three ensembles represent the same induced density $\mu[v]$ and differ only in the location distributions where the particles are placed. All ensembles start at $t=0$ and are advanced with the trained velocity. As for the training data, in window $I_j$, only the ensemble corresponding to the period that contains the midpoint of $I_j$ is used for training. We call it the training ensemble of $I_j$. Each window is therefore trained on particles placed where the mass of its own period comes from. We propagate the other two ensembles through the window $I_j$ using the trained velocity. More details can be found in Appendix~\ref{app:banks}.

The adaptive weights below address the following problem: a region that holds little mass enters $J_{\mathrm{reg}}$ with that small mass, whereas the KL divergence it produces is later multiplied by the growth of the region. The adaptive weights give such particles extra weight before their own mass has grown.

Let $\widehat\lw_i(t)$ denote the log-weight of particle $i$ along the flow of $\rho_{\mathrm{pred}}$ (Appendix~\ref{app:banks}), and $\Ghat_i^{(j)}:=\widehat\lw_i(T)-\widehat\lw_i(T_j)$ its predicted growth from the start of window $I_j=[T_j,T_{j+1}]$ to the endpoint $T$. In both terms of \eqref{eq:J_modular}, we replace the physical weights by the adaptive weights
\begin{equation}\label{eq:measure}
  \meas^i_n=\frac12\,\frac{e^{\lw^i_n}}{\sum_{i'=1}^Ne^{\lw^{i'}_n}}
            +\frac12\,\frac{e^{\lw^i_n+\Ghat_i^{(j)}}}{\sum_{i'=1}^Ne^{\lw^{i'}_n+\Ghat_{i'}^{(j)}}},
\end{equation}
where the sums run over the $N$ particles of the training ensemble. The first term weights a particle by its current mass and the second by its current mass multiplied by its predicted growth from $T_j$ to $T$. Since every particle keeps at least half of its physical weight, the loss with physical weights is at most twice the loss with adaptive weights.
\subsection{Windowed iteration and the complete algorithm}\label{sec:solver}
For every window $I_j$, we will train a corresponding neural network $v_{\theta_j}$. In the actual algorithm, we train the time windows sequentially. After training the previous window, we warm-start the next window from the velocity $v_{\bar{\theta}^{\mathrm{final}}}$ learned on the previous window and use the particle states at the end of the previous window as the initial conditions of the next window. The particle states at $T_j$ are therefore produced by velocities already trained on the earlier windows, so the iteration in each window only has to resolve the dynamics over a short time interval.

With physical weights and a quadratic loss, the regression term of the frozen iteration of Section~\ref{sec:iter} has the same structure as the loss of~\cite{li2023scvm}, which takes one gradient step per iteration. Here the targets contain the carried score, the characteristic system contains the reaction, and the iteration runs window by window with the stopping rule of Section~\ref{sec:iter}.

\begin{algorithm}[t]
\caption{Windowed self-consistent training with frozen velocities}\label{alg:main}

\begin{algorithmic}[1]
\Require $\rho_0,b,\gamma,g,T$; windows $I_j=[T_j,T_{j+1}]$, $j=0,\dots,N_w-1$; period-adapted sampling distributions $q_0^m$ for the periods $\mathcal T_m$, $m=0,1,2$; regularization weight $\lamres$
\Ensure window velocities $v_{\theta_0},\dots,v_{\theta_{N_w-1}}$
\For{$m=0,1,2$} \Comment{one ensemble per period}
  \State draw $X_0^{i}\sim q_0^m$; set $\ell_0^{i}=\log\rho_0(X_0^{i})$, $\Sc_0^{i}=\nabla\log\rho_0(X_0^{i})$, $\lw_0^{i}=\log\rho_0(X_0^{i})-\log q_0^m(X_0^{i})$
\EndFor
\State pretrain $v_\theta$ on the exact velocity $b-\gamma\nabla\log\rho_0$ at $t=0$
\For{$j=0,\dots,N_w-1$}
  \State select the training ensemble, whose period contains the midpoint of $I_j$
  \Repeat \Comment{outer iterations}
    \State freeze the velocity, $\bar\theta\gets\theta$
    \State integrate the training ensemble over $I_j$ with $v_{\bar\theta}$ by \eqref{eq:heun}, from its state at $T_j$
    \State labels $v_{\rm label}=b-\gamma\Sc$ and adaptive weights $\meas$ from $\lw$ and $\Ghat$ by \eqref{eq:measure}
    \State update $\theta$ by a fixed number of optimizer steps (Section~\ref{sec:frozen})
  \Until{$v_\theta\approx v_{\bar\theta}$ (Appendix~\ref{app:iteration}), or after $k_{\max}$ outer iterations}
  \State set $\theta_j\gets\theta$ and advance every ensemble over $I_j$ with $v_{\theta_j}$
\EndFor
\end{algorithmic}
\end{algorithm}

\section{Numerical experiments}\label{sec:numerics}

The experiments are organized around three types of equations: the OU equation with a linear source, that is, the Fokker--Planck equation of an OU process with an added linear reaction term, the Fisher--KPP equation, and WFR gradient flows of the KL divergence. The OU case tests the accuracy of the carried density, score and weights. The Fisher--KPP equation has a reaction term that depends on the density, so errors in the carried density will feed back to the system. In the WFR gradient flow, the target can have modes far from the initial mass, and these modes must grow from regions that initially carry almost no mass.

We use Algorithm~\ref{alg:main} with shared network and integration settings, listed in Table~\ref{tab:settings} of Appendix~\ref{app:settings}. The particles are integrated by Heun's method with $\Delta t=0.0025$ in windows of width $W=0.01$. The regularizer term has the weight $\lamres=0.1$. We report the relative velocity error $E_v$ over the report times, the relative mass error $e_M$, the distribution error $D$, and the relative carried-score error $E_S$ against independent reference solutions. Appendix~\ref{app:evaluation} defines these quantities.

\subsection{OU equation with a linear source}\label{sec:ou}

We first validate the method on a problem with an analytical solution.
Consider the drift $b(x)=-\kappa x$ with the linear source $g(x)=a^\top x$,
\begin{equation}\label{eq:ou_pde}
  \partial_t\rho=\nabla\cdot(\kappa x\,\rho)+\gamma\Delta\rho+(a^\top x)\,\rho .
\end{equation}

For $\kappa,\gamma>0$ and initial data $\rho_0=M_0\,\mathcal N(m_0,C_0)$ with $C_0$ positive definite, the exact solution is
\begin{equation}\label{eq:ou_density}
  \rho(t,x)
  =\frac{M(t)}{(2\pi)^{d/2}\sqrt{\det C(t)}}
   \exp\!\left[-\frac12\bigl(x-m(t)\bigr)^\top
   C(t)^{-1}\bigl(x-m(t)\bigr)\right].
\end{equation}
Writing $r_t=e^{-\kappa t}$, its covariance, mean and total mass are given by
\begin{equation}\label{eq:ou_exact}
  \begin{aligned}
    C(t)
      &=r_t^2C_0+\frac{\gamma}{\kappa}(1-r_t^2)I,\\[3pt]
    m(t)
      &=r_t m_0+\frac{r_t(1-r_t)}{\kappa}\,C_0a
        +\frac{\gamma}{\kappa^2}(1-r_t)^2a,\\[3pt]
    \log\frac{M(t)}{M_0}
      &=\frac{1-r_t}{\kappa}\,a^\top m_0
        +\frac{(1-r_t)^2}{2\kappa^2}\,a^\top C_0a
        +\frac{\gamma}{2\kappa^3}
         (2\kappa t-3+4r_t-r_t^2)|a|^2.
  \end{aligned}
\end{equation}
The corresponding score and velocity are
\begin{equation}\label{eq:ou_fields}
  \begin{aligned}
    \nabla\log\rho(t,x)&=-C(t)^{-1}\bigl(x-m(t)\bigr),\\
    v_\rho(t,x)&=-\kappa x+\gamma C(t)^{-1}\bigl(x-m(t)\bigr).
  \end{aligned}
\end{equation}
The linear source changes the mean and total mass, and covariance is not affected by the source term. Appendix~\ref{app:ou_exact} derives these formulas.

We take $d=2$, $\kappa=\gamma=1$, $a=(-1,-1)^\top$, $M_0=1$, $m_0=(1,-0.5)^\top$, $C_0=\bigl(\begin{smallmatrix}1.4&0.25\\0.25&0.8\end{smallmatrix}\bigr)$ and $T=1$. The report times are $0.02$, $0.05$, $0.1$, $0.2$, $0.4$, $0.6$, $0.8$ and $1$, and the distribution error uses $16\times16$ grids on $[-25,25]^2$.

\begin{figure}[tbp]
\centering
\includegraphics[width=\textwidth]{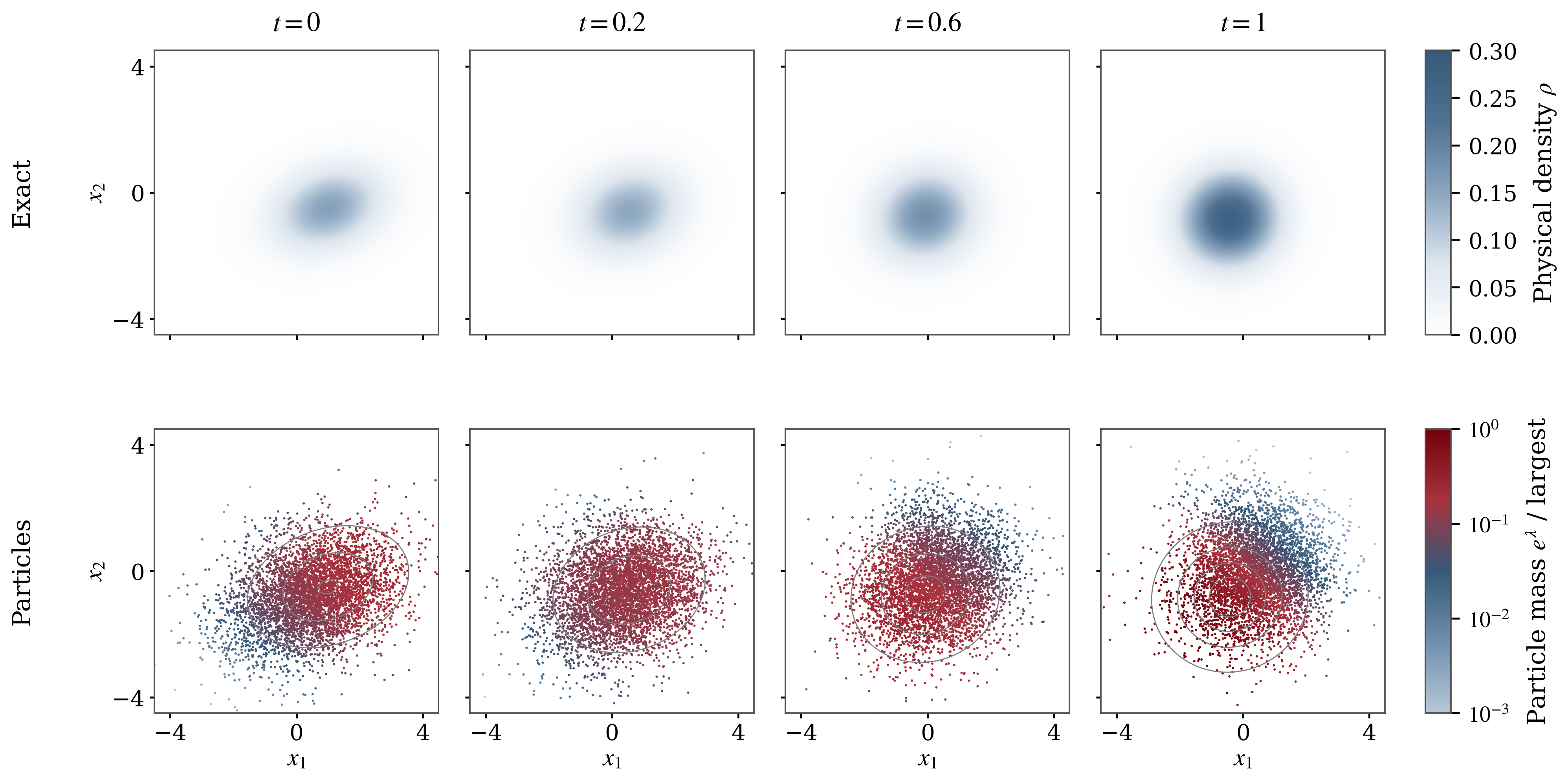}
\caption{OU density snapshots at $t=0$, $0.2$, $0.6$ and $1$. Top: exact density. Bottom: 5000 output particles, colored by the mass $e^{\lambda}$ each one carries (log scale, relative to the heaviest). Particle positions follow the sampling distribution; the density is the mass-weighted particle measure (Appendix~\ref{app:metrics}).}
\label{fig:ou_snapshots}
\end{figure}

\begin{figure}[tbp]
\centering
\includegraphics[width=\textwidth]{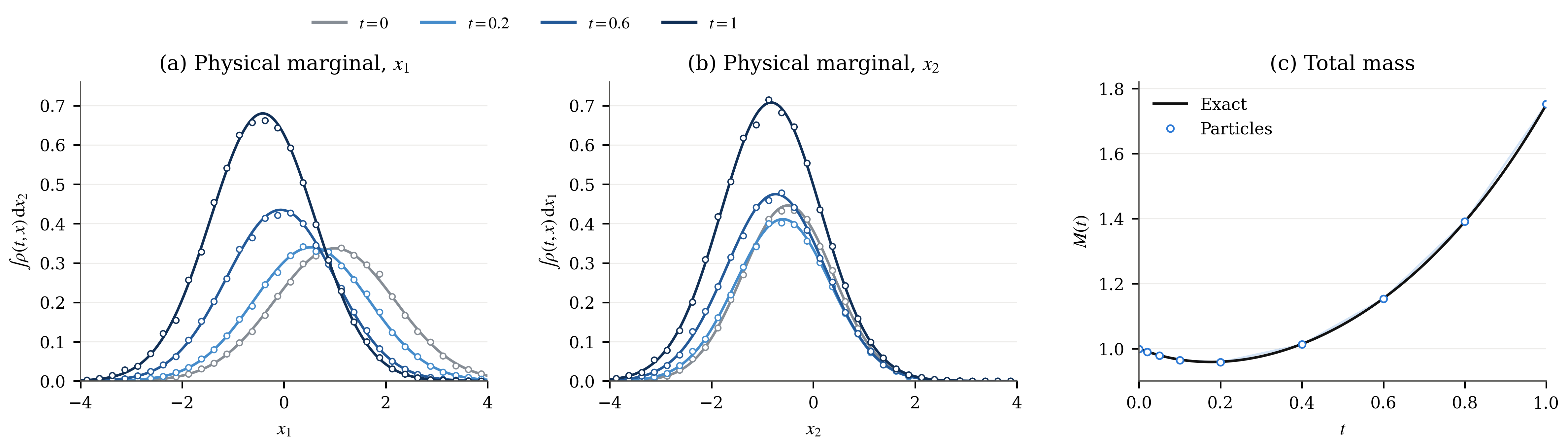}
\caption{OU equation with a linear source. (a,b) Physical coordinate marginals at $t=0$, $0.2$, $0.6$ and $1$: exact marginals (lines) and weighted particle histograms (dots) of the output particles. (c) Exact mass (line) and the median particle estimate (dots).}
\label{fig:ou_solution}
\end{figure}

Figure~\ref{fig:ou_snapshots} shows the two-dimensional density at successive times. Figure~\ref{fig:ou_solution} compares the two coordinate marginals with those of the exact solution at four time points and shows the change of the total mass. The mass error, which tests the carried weights, and the carried-score error stay below $1\%$, and the distribution error below $0.01$ (Table~\ref{tab:smooth}). The carried score is driven by the derivatives $\nabla v_\theta$ and $\nabla(\nabla\cdot v_\theta)$ of the learned field (Proposition~\ref{prop:carry}), and their root-mean-square errors are at most $1.6\times10^{-2}$ at $t=0$ and $4\times10^{-4}$ at $t=1$ (Appendix~\ref{app:metrics}).

\FloatBarrier

\subsection{Fisher--KPP equation}\label{sec:kpp}

We next consider the Fisher--KPP equation
\begin{equation}\label{eq:kpp_pde}
  \partial_t\rho=\gamma\Delta\rho+k\,\rho\Bigl(1-\frac{\rho}{K}\Bigr).
\end{equation}
The reaction term explicitly depends on the density: it is positive where $\rho<K$ and negative where $\rho>K$, and its rate is evaluated at the carried density $e^{\ell}$. In a biological context, $K$ means the carrying capacity of the environment. According to \eqref{eq:char_lw} and \eqref{eq:nabla_g}, the error in the carried log-density feeds back into both the weight and the score. We consider two experiments with $\gamma=1$, $K=1$ and $T=2$:
\begin{enumerate}
  \item \textbf{One-dimensional problem}: $\rho_0(x)=4\,\mathcal N(0,1)=\frac{4}{\sqrt{2\pi}}e^{-\frac{x^2}{2}}$ and $k=1$.
  \item \textbf{Three-dimensional problem}: $\rho_0(x)=\frac12\exp\bigl(-\frac{|x|^2}{2\cdot0.8^2}\bigr)$ on $\mathbb R^3$, $k=2$.
\end{enumerate}

\begin{figure}[t]
\centering
\includegraphics[width=\textwidth]{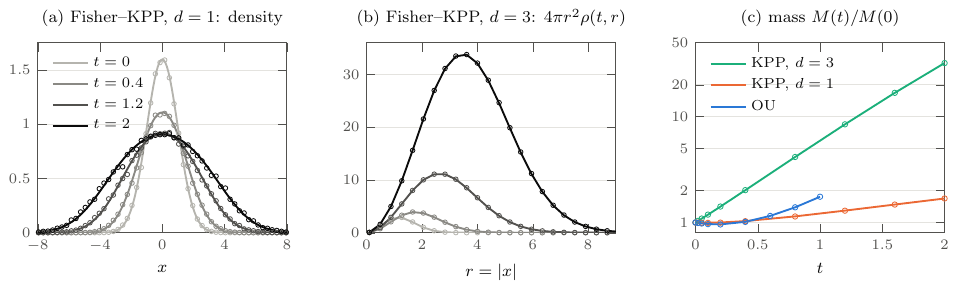}
\caption{Fisher--KPP profiles and mass evolution. (a) Fisher--KPP in one dimension: reference density (lines) and physically weighted particle histogram (circles) at four time points. (b) Fisher--KPP in three dimensions: radial mass density $4\pi r^2\rho$ of the reference (lines) and of the carried density averaged over radial bins. (c) Total mass relative to the initial mass for the two Fisher--KPP problems and the OU equation of Section~\ref{sec:ou}: reference (lines) and particle estimate $\widehat M$ (circles).}
\label{fig:accuracy}
\end{figure}

The two initial conditions show different parts of the Fisher--KPP dynamics. In one dimension, the initial peak $4/\sqrt{2\pi}\approx1.6$ exceeds the capacity $K$, so the reaction term will decrease mass at the center while the tail region will grow where $\rho<K$.
The peak decreases, and at $t=2$ the distribution is a broad peak slightly below the capacity whose flanks spread outwards (Figure~\ref{fig:accuracy}a); the mass grows from $4$ to $6.74$ (Figure~\ref{fig:accuracy}c).
In three dimensions the initial density $\rho_0$ is below the capacity everywhere, thus the solution grows and spreads at the same time. The total mass grows from $4.0$ to $129$ (Figure~\ref{fig:accuracy}c), and the maximum of the radial mass density $4\pi r^2\rho$ moves outwards (Figure~\ref{fig:accuracy}b).

In both problems the particle solutions match the reference solutions over the whole time interval, with relative velocity, mass and carried-score errors below $0.3\%$ (Table~\ref{tab:smooth} and Figure~\ref{fig:accuracy}). In one dimension the weighted histogram reproduces the decay of the over-capacity peak and the outward spreading of the flanks; in three dimensions the carried density reproduces both the bulk and the outgoing front along the radius.

\begin{figure}[tbp]
\centering
\includegraphics[width=\textwidth]{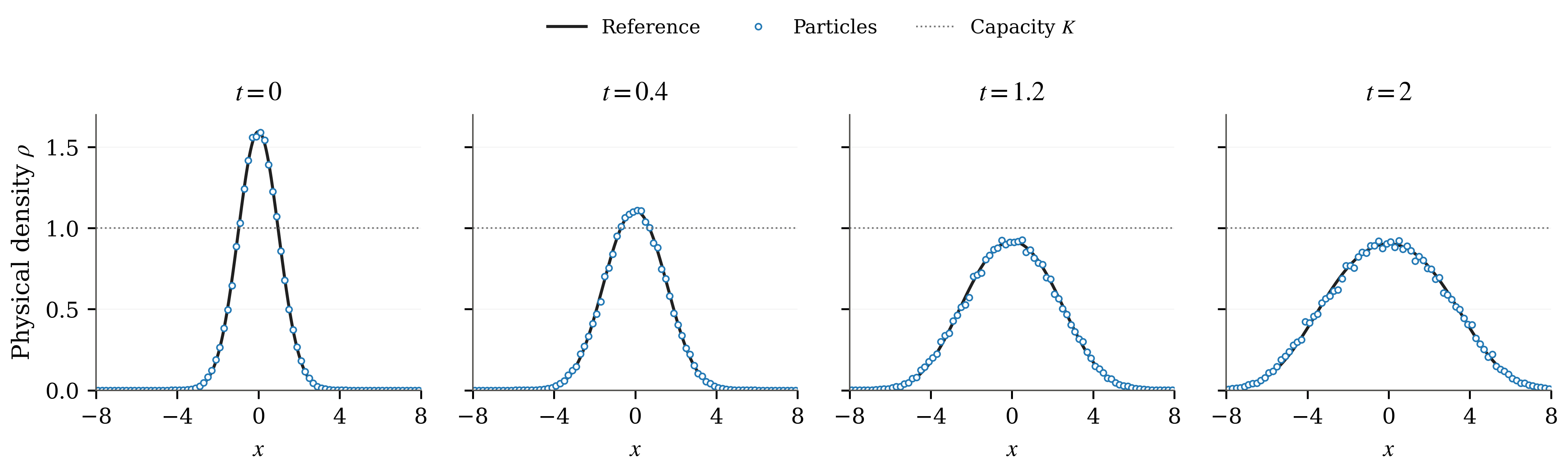}
\caption{One-dimensional Fisher--KPP snapshots. Each panel compares the reference density (line) and the physical weighted histogram at one time. The dotted line marks the carrying capacity $K=1$.}
\label{fig:kpp1d_snapshots}
\end{figure}

\begin{figure}[tbp]
\centering
\includegraphics[width=\textwidth]{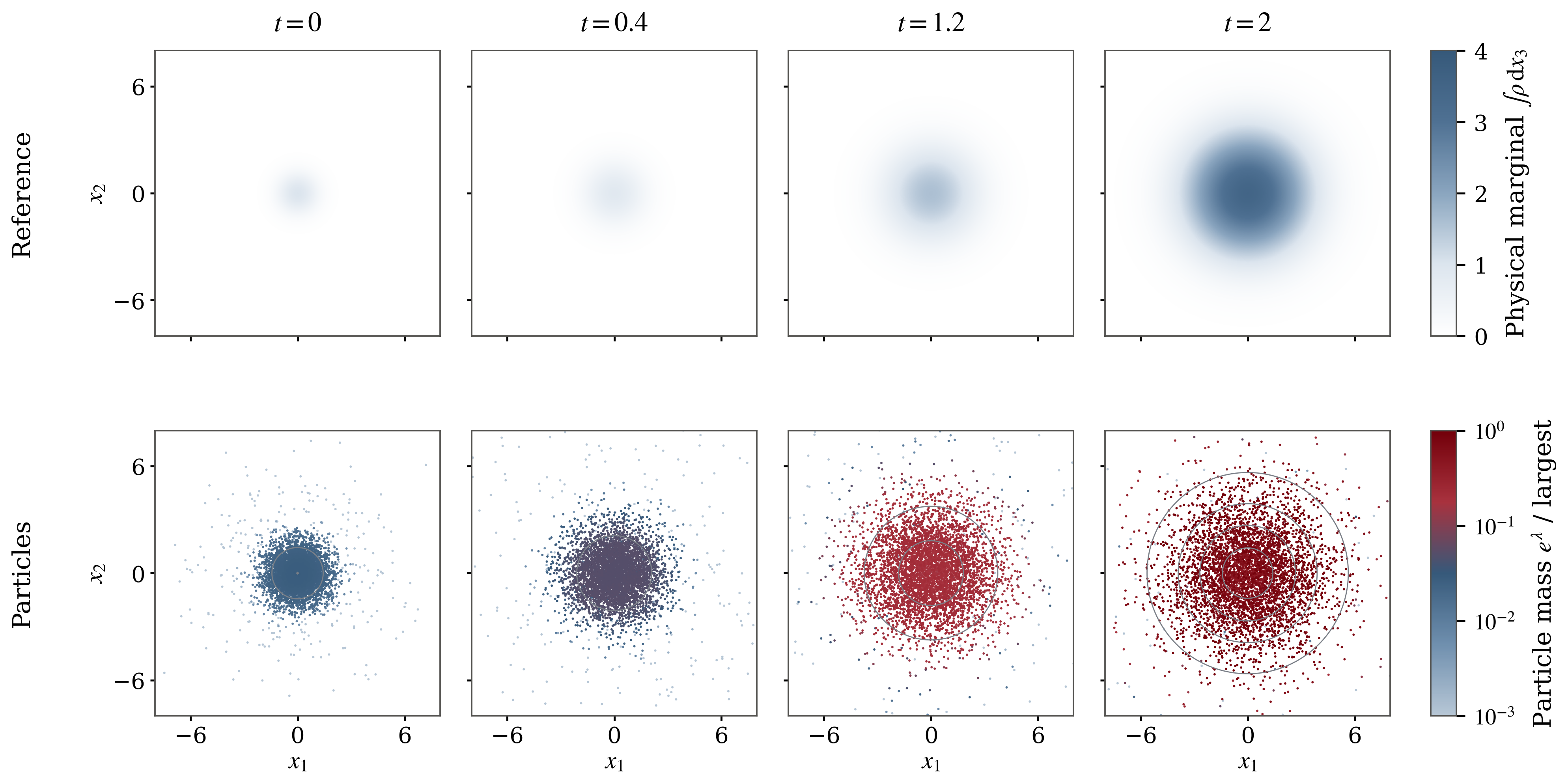}
\caption{Three-dimensional Fisher--KPP snapshots of $\int\rho(t,x_1,x_2,x_3)\,\mathrm dx_3$. Top: the radial reference integrated in $x_3$. Bottom: 5000 of the three-dimensional output particles projected onto $(x_1,x_2)$, colored by their carried mass $e^{\lambda}$ (log scale, relative to the heaviest); the color change records the growth of mass along the characteristics.}
\label{fig:kpp3d_snapshots}
\end{figure}

Figures~\ref{fig:kpp1d_snapshots} and~\ref{fig:kpp3d_snapshots} separate the evolution into individual times. The three-dimensional panels show the physical $(x_1,x_2)$ marginal obtained by integrating out $x_3$.

\begin{table}[t]
\centering
\caption{OU and Fisher--KPP: errors against reference solutions. $E_v$: range over the training seeds; $e_M$, $D$ and $E_S$: maximum values over the report times, training seeds and the two output seeds (Appendix~\ref{app:metrics}).}
\label{tab:smooth}
\smallskip
\small
\begin{tabular}{@{}lccccc@{}}
\toprule
Problem & $d$ & $E_v$ (\%) & $e_M$ (\%) & $D$ & $E_S$ (\%) \\
\midrule
OU with linear source & 2 & 0.145--0.152 & 0.56 & 0.0094 & 0.05 \\
Fisher--KPP & 1 & 0.10--0.18 & 0.27 & 0.026 & 0.11 \\
Fisher--KPP & 3 & 0.11--0.12 & 0.17 & 0.030 & 0.08 \\
\bottomrule
\end{tabular}
\end{table}

\FloatBarrier

\subsection{Wasserstein--Fisher--Rao gradient flows of the KL divergence}\label{sec:wfr}

We finally solve the Wasserstein--Fisher--Rao gradient flow of the KL divergence towards a target density distribution $\pi$,
\begin{equation}\label{eq:klwfr_pde}
  \partial_t\rho=\gamma\nabla\cdot\Bigl(\rho\,\nabla\log\frac{\rho}{\pi}\Bigr)-\beta\,\rho\log\frac{\rho}{\pi},
\end{equation}
which is \eqref{eq:rd_general} with $b=\gamma\nabla\log\pi$ and $g(x,\rho)=-\beta\log(\rho/\pi(x))$.
The transport part moves mass along $-\nabla\log(\rho/\pi)$, while the reaction part increases mass when $\rho<\pi$ and decreases it when $\rho>\pi$. The total mass is not conserved: $\frac{\mathrm d}{\mathrm dt}M=-\beta\int\rho\log(\rho/\pi)\dd x$.
Since $\rho$ is not normalized, \eqref{eq:klwfr_pde} is the WFR gradient flow of the generalized KL divergence $\int\bigl(\rho\log(\rho/\pi)-\rho+\pi\bigr)\dd x$ on nonnegative densities. Here $\rho_0$ and $\pi$ have unit mass, $\gamma=1$, $\beta=4$ and $T=1$.

\subsubsection*{Mixed Gaussian target}\label{sec:wide}

The initial density is a single Gaussian, $\rho_0=\mathcal N\bigl((0,8)^\top,\,0.3I\bigr)$, and the target is a mixture of four Gaussians with equal weights,
\begin{gather*}
  \pi=\tfrac14\sum_{k=1}^4\mathcal N(m_k,\Sigma_k),\qquad
  m_1=(0,8)^\top,\quad m_2=(0,2)^\top,\quad m_{3}=(-3,5)^\top,\quad m_4=(3,5)^\top,\\
  \Sigma_1=\Sigma_2=\operatorname{diag}(1.2,\,0.3),\qquad \Sigma_3=\Sigma_4=\operatorname{diag}(0.3,\,2).
\end{gather*}
The report times are $0.02$, $0.05$, $0.1$, $0.15$, $0.2$, $0.25$, $0.3$, $0.4$, $0.5$, $0.6$, $0.8$ and $1$.
At the bottom component, $\rho_0$ is far below $\pi$, so $\log(\rho/\pi)$ is strongly negative and the reaction is large there. The bottom mode is formed from the far tail of $\rho_0$, whose weight the reaction amplifies by many orders of magnitude: its share of the total mass is $4\times10^{-4}$ at $t=0.3$, $2\%$ at $t=0.5$ and $20\%$ at $t=1$. Meanwhile the total mass first falls from $1$ to $0.38$ at $t\approx0.4$, because the reaction term removes much of the mass in the top mode.

The particles that end in the bottom mode carry almost no weight at the beginning. The period-adapted sampling distributions place initial particles where this later mass comes from, and the adaptive weights give them weight in the regression before their own mass has grown (Section~\ref{sec:banks}).

\begin{figure}[tbp]
\centering
\includegraphics[width=0.8\textwidth]{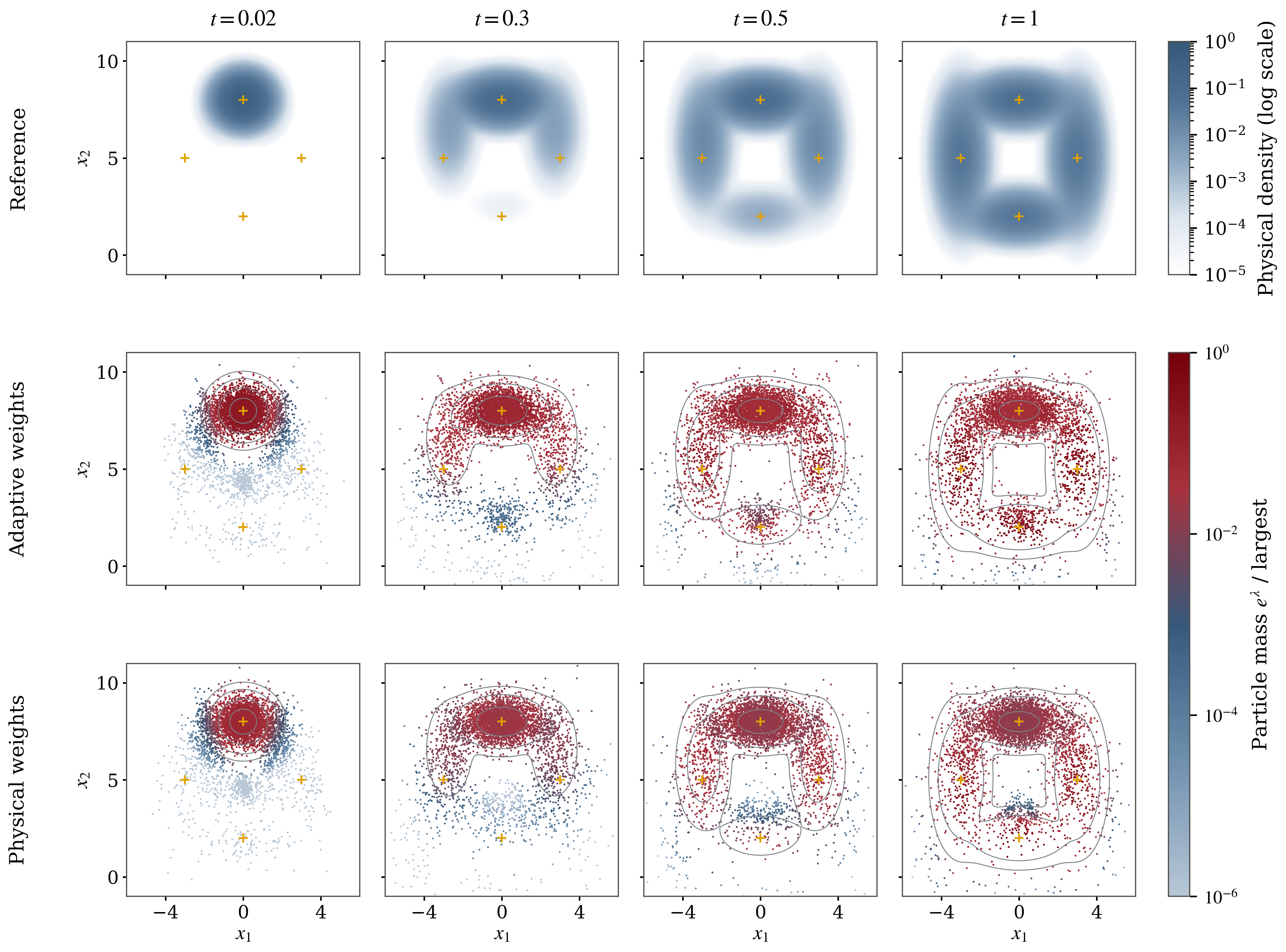}
\caption{Two-dimensional snapshots for the mixed Gaussian target. Top: reference density (log scale). Middle and bottom: 5000 output particles trained with adaptive and physical regression weights, colored by the mass $e^{\lambda}$ each one carries (log scale, relative to the heaviest). With adaptive weights the particles that will form the bottom mode are already present at $t=0.02$, with almost no mass, and gain it by $t=0.5$. Crosses mark the centers of the four target modes.}
\label{fig:wide_snapshots}
\end{figure}

\begin{figure}[tbp]
\centering
\includegraphics[width=\textwidth]{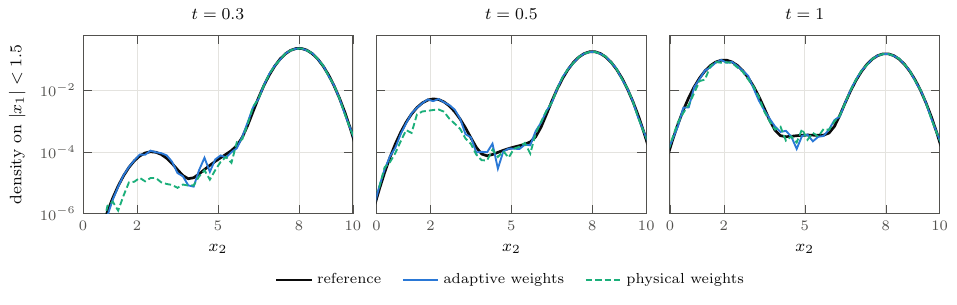}
\caption{Mixed Gaussian target: density on the central strip $|x_1|<1.5$ as a function of $x_2$, for the reference and for the particle solutions trained with adaptive and with physical weights (one training seed). The peak near $x_2=2$ is the bottom mode.}
\label{fig:wide}
\end{figure}

\begin{figure}[t]
\centering
\includegraphics[width=0.75\textwidth]{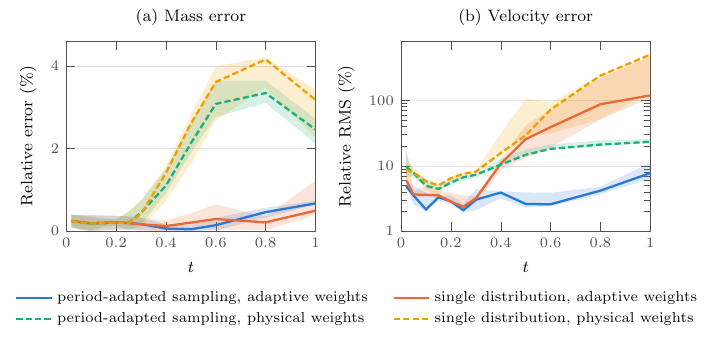}
\caption{Mixed Gaussian target: initial sampling distributions and regression weights. (a) Relative mass error; lines are medians and bands ranges over the three training seeds and the two output seeds. (b) Velocity error at each report time, divided by the root-mean-square of the exact velocity at that time; lines are medians and bands ranges over the three training seeds. The exact velocity decreases by a factor of about 25 from $t=0.02$ to $t=1$, so the late times contribute little to the pooled error $E_v$ of Table~\ref{tab:factorial}.}
\label{fig:wide_factorial}
\end{figure}

Figure~\ref{fig:wide_snapshots} (middle row) shows these particles: at $t=0.02$ they already lie below the top mode with almost no mass, and they gain it as the reference mode forms, by $t=0.5$. On the central strip the particle density follows the reference over several decades, including the small bottom peak at $t=0.3$ (Figure~\ref{fig:wide}). Over three training seeds the velocity error $E_v$ is about $4\%$ and the mass error at most $0.74\%$ (Table~\ref{tab:factorial} and Figure~\ref{fig:wide_factorial}, blue); the carried-score error is at most $1.3\%$.

\paragraph{Comparison of sampling distributions and regression weights.}
We design a comparison experiment to test the period-adapted sampling and adaptive weights mechanisms. The particles are drawn either from the three period-adapted sampling distributions or from a single distribution used in every window (Appendix~\ref{app:banks}). The regression weights are either the adaptive weights \eqref{eq:measure} or the physical weights \eqref{eq:basic_choice}. All other settings and the training seeds are the same, and all runs are evaluated on the same output particles.

With the physical weights, the bottom mode forms later than in the reference. The mass error reaches $3.6$--$4.2\%$ and the distribution error $D$ is $0.049$--$0.063$, against at most $1.2\%$ and $0.014$ with adaptive weights (Table~\ref{tab:factorial} and Figure~\ref{fig:wide_factorial}). In every training seed, the adaptive weights give the smaller mass error and distribution error, with either sampling distribution.
The sampling distributions act mainly on the velocity field: the period-adapted ones give the smaller velocity error.

\begin{table}[t]
\centering
\caption{Mixed Gaussian target: initial sampling distribution and regression weights. $E_v$: range over three training seeds; $e_M$ and $D$: worst values over the report times, training seeds and the two output seeds.}
\label{tab:factorial}
\smallskip
\small
\begin{tabular}{@{}llccc@{}}
\toprule
Sampling distribution & Weights & $E_v$ (\%) & $e_M$ (\%) & $D$ \\
\midrule
period-adapted & adaptive & 3.6--4.2 & 0.74 & 0.013 \\
single & adaptive & 6.3--20.4 & 1.2 & 0.014 \\
period-adapted & physical & 7.9--10.3 & 3.6 & 0.049 \\
single & physical & 10.8--25.6 & 4.2 & 0.063 \\
\bottomrule
\end{tabular}
\end{table}

\FloatBarrier
\section{Conclusion and discussion}\label{sec:conclusion}

We developed a score-based neural ODE method for solving reaction-diffusion-advection equations with a local multiplicative reaction term. The particles carry their log-density, score and physical weights along trajectories rather than estimating them. A relative-entropy inequality motivates the mass-normalized self-consistency residual which is used as the continuous objective. The numerical algorithm learns a velocity field through successive regressions on particle data, with period-adapted initial sampling distributions and adaptive regression weights.

The reaction changes the local density and the relative importance of different regions. Therefore, in a deterministic particle method, the density, the score and the physical weights must all be available at the particles. We propagate them along the trajectories by ODEs. Estimating the score from the particle ensembles, by kernels~\cite{maoutsa2020interacting,zhou2024forwardbackward} or by local regression~\cite{willems2026sformulation}, adds an estimation error that depends on the number and the arrangement of the particles. When the score is instead obtained by integrating its ODE along each trajectory (Proposition~\ref{prop:carry}), the error will also depend on the accuracy of the time integration. The computation cost is also higher because the spatial derivatives of the velocity must be computed on every particle. Moreover, the carried score is the score of the density induced by the current velocity, so it approaches the score of the solution only as the velocity approaches self-consistency. In training, each outer iteration is a regression on labels from a frozen velocity. These labels lag behind the current velocity, so several iterations are needed in each window, which increases the computational cost.

A region that becomes important later may contain only a few particles in the beginning. Sequential Monte Carlo methods handle this by resampling, branching or killing particles~\cite{delmoral2004feynman,delmoral2006smc}. We choose to change the initial sampling distribution instead, which leaves the induced density $\mu[v]$ unchanged (Lemma~\ref{lem:weighted_residual}). Every particle is weighted by its physical mass in the objective functional, so regions with small physical mass receive less weight in training. We therefore modify the regression weights so that particles in regions that become important later can have more influence in the training. The period-adapted sampling distributions are designed to cover the regions of interest in each period, while the adaptive weights make the optimizer prioritize them.

Numerical results demonstrate the method on problems with changing mass and growing modes. The WFR gradient flow example shows why the allocation of training effort matters: particles with little physical weight at early time can contribute substantially later. The adaptive weights and the period-adapted sampling distributions can improve the accuracy of the solution (Table~\ref{tab:factorial} and Figure~\ref{fig:wide_factorial}).

The method has several limitations. The error bound of Theorem~\ref{thm:stability}(ii) replaces the growth rate by its largest value. For the OU process and the WFR gradient flow the growth rate is unbounded, so part (ii) does not apply and the accuracy there is assessed only numerically. Moreover, no numerical analysis yet connects the scheme to the bound. Our experiments do not include singularly perturbed problems with a moving thin layer; for such problems the particles have to resolve the layer, and the network has to fit the large high-order derivatives that enter the score equation. Finally, the integration cost grows as the dimension increases.

Several directions remain open. The characteristics could be combined with particle methods that add or remove particles adaptively. The sampling distributions and the adaptive weights could be built from a short preliminary run of a sampling method instead of a problem-specific predicted density. We will study the selection of temporal integrators for the proposed neural ODE systems. On the theory side, sharper bounds and the bridge between the continuous objective and the trained loss are left for future work.

\section*{Declaration of generative AI use}
During this work, the authors used OpenAI's ChatGPT and Codex and Anthropic's Claude to assist with language and \LaTeX{} editing and with the implementation and debugging of the numerical code. These tools also assisted with the organization and presentation of the computational results. The authors reviewed, verified and revised all AI-assisted material and take full responsibility for the mathematical statements and proofs, the numerical experiments and the conclusions of this paper.

\section*{Acknowledgments}
{\sloppy Jianguo Hou and Wuchen Li are supported in part by the AFOSR YIP award No.\ FA9550-23-1-0087, NSF RTG:~2038080, NSF FRG:~2245097, NSF DMS:~2601952, and the McCausland Faculty Fellowship at the University of South Carolina.\par}

\bibliographystyle{plain}
\bibliography{refs}

\appendix

\section{Characteristics and weighted particles}\label{app:proofs}

This appendix derives the characteristic system of Proposition~\ref{prop:carry} and the weighted identity of Lemma~\ref{lem:weighted_residual}, and proves the two propositions on the self-consistency residual. Throughout, $v$ is an admissible velocity field (Definition~\ref{def:admissible}) and Assumption~\ref{ass:continuum} holds. All statements concern strictly positive classical solutions, that is, solutions that are $C^2$ in $x$, continuous in $t$ and $C^1$ in $t$ on each window.

\subsection{The flow and the log-density}\label{app:proof_carry}

On each window of Definition~\ref{def:admissible}, $v$ is continuous in $t$ and $\nabla v$ is bounded, so $v$ is Lipschitz in $x$ uniformly in $t$. Solving \eqref{eq:char_x} window by window therefore gives a unique trajectory $X_t=\Phi_t(x_0)$ on $[0,T]$ for every starting point $x_0\in\R^d$, and the flow map $\Phi_t$ is a $C^3$ diffeomorphism of $\R^d$. Its Jacobian determinant $J_t(x_0)=\det\nabla\Phi_t(x_0)$ satisfies Liouville's formula
\begin{equation}\label{eq:app_liouville}
  \frac{\mathrm d}{\mathrm dt}\log J_t(x_0)=\nabla\!\cdot v\bigl(t,\Phi_t(x_0)\bigr),\qquad J_0(x_0)=1,
\end{equation}
on each window, and $J_t$ is continuous across the window interfaces, so $J_t>0$.

Let $\mu>0$ be a classical solution of the transport-reaction equation \eqref{eq:mu_pde}. Expanding the divergence,
\[
  \partial_t\mu=-\nabla\!\cdot(\mu v)+g(t,x,\mu)\,\mu
  =-\mu\,\nabla\!\cdot v-v\cdot\nabla\mu+g(t,x,\mu)\,\mu ,
\]
and dividing by $\mu$ gives
\begin{equation}\label{eq:app_logmu}
  \partial_t\log\mu=-\nabla\!\cdot v-v\cdot\nabla\log\mu+g(t,x,\mu).
\end{equation}
Along a trajectory, the chain rule gives
\[
  \frac{\mathrm d}{\mathrm dt}\log\mu(t,X_t)=\partial_t\log\mu(t,X_t)+v(t,X_t)\cdot\nabla\log\mu(t,X_t).
\]
Substituting \eqref{eq:app_logmu}, the transport term $v\cdot\nabla\log\mu$ cancels, and we obtain:
\[
  \frac{\mathrm d}{\mathrm dt}\log\mu(t,X_t)=-\nabla\!\cdot v(t,X_t)+g\bigl(t,X_t,\mu(t,X_t)\bigr),
\]
so $\log\mu(t,X_t)$ solves \eqref{eq:char_ell}.

Conversely, we construct $\mu$ from \eqref{eq:char_ell}. Fix $x_0$ and write
\[
  F_{x_0}(t,\ell)=-\nabla\!\cdot v\bigl(t,\Phi_t(x_0)\bigr)+g\bigl(t,\Phi_t(x_0),e^{\ell}\bigr).
\]
By Assumption~\ref{ass:continuum}, $F_{x_0}$ is continuous in $t$ on each window and locally Lipschitz and non-increasing in $\ell$, so \eqref{eq:char_ell} with $\ell_0=\log\rho_0(x_0)$ has a unique local solution. Monotonicity prevents blow-up: while $\ell_t\ge\ell_0$ we have $\dot\ell_t\le F_{x_0}(t,\ell_0)$, and while $\ell_t\le\ell_0$ we have $\dot\ell_t\ge F_{x_0}(t,\ell_0)$, so
\[
  |\ell_t-\ell_0|\le\int_0^t|F_{x_0}(s,\ell_0)|\dd s ,
\]
and the solution exists on $[0,T]$. Define
\[
  \mu(t,x):=\exp\bigl(\ell_t(\Phi_t^{-1}(x))\bigr).
\]
Since $v$ is $C^3$ in $x$ and $g$ is $C^2$ in $(x,r)$, $\mu$ is strictly positive, $C^2$ in $x$, continuous in $t$ and $C^1$ in $t$ on each window. Reading the computation above backwards, \eqref{eq:char_ell} says that $\mu$ solves \eqref{eq:mu_pde} on each window, and $\mu(0)=\rho_0$ because $\ell_0=\log\rho_0(X_0)$. The solution is unique among strictly positive classical solutions: for another one, $\tilde\mu$, the function $\log\tilde\mu(t,\Phi_t(x_0))$ solves \eqref{eq:char_ell} with the same initial value, so $\tilde\mu=\mu$. In particular $\ell_t=\log\mu(t,X_t)$ along every trajectory, which is the first half of Proposition~\ref{prop:carry}(i). The weight $\lw_t$ then follows from \eqref{eq:char_lw} by integration.

\subsection{Score dynamics}\label{app:score}

Write $\ell(t,x)=\log\mu(t,x)$ for the log-density as a field, so that \eqref{eq:app_logmu} reads
\begin{equation}\label{eq:app_ell_field}
  \partial_t\ell+v\cdot\nabla\ell=-\nabla\!\cdot v+g\bigl(t,x,e^{\ell}\bigr).
\end{equation}
We take the gradient of each term. Since $\nabla^2\ell$ is symmetric,
\begin{align*}
  \nabla\bigl(v\cdot\nabla\ell\bigr)&=(\nabla v)^\top\nabla\ell+(\nabla^2\ell)\,v,\\
  \nabla\bigl(-\nabla\!\cdot v\bigr)&=-\nabla(\nabla\!\cdot v),\\
  \nabla\bigl[g\bigl(t,x,e^{\ell}\bigr)\bigr]&=\partial_xg+\partial_\rho g\,e^{\ell}\,\nabla\ell .
\end{align*}
The gradient of \eqref{eq:app_ell_field} is therefore
\[
  \partial_t\nabla\ell+(\nabla^2\ell)\,v
  =-(\nabla v)^\top\nabla\ell-\nabla(\nabla\!\cdot v)+\partial_xg+\partial_\rho g\,e^{\ell}\,\nabla\ell .
\]
Along a trajectory,
\[
  \frac{\mathrm d}{\mathrm dt}\nabla\ell(t,X_t)=\partial_t\nabla\ell+(\nabla^2\ell)\,v ,
\]
which is exactly the left-hand side. Hence $\nabla\ell(t,X_t)$ solves the score equation \eqref{eq:char_score} with the total gradient \eqref{eq:nabla_g}, and its initial value is $\nabla\log\rho_0(X_0)=\Sc_0$. Equation \eqref{eq:char_score} is linear in $\Sc$ with coefficients that are continuous on each window, so its solution is unique, and
\[
  \Sc_t=\nabla\log\mu(t,X_t).
\]
This completes Proposition~\ref{prop:carry}(i). Since the term $(\nabla^2\ell)\,v$ is absorbed by the motion of the particle, \eqref{eq:char_score} contains no second derivatives of $\ell$, and the carried score requires only $\nabla v$ and $\nabla(\nabla\!\cdot v)$.

\subsection{Weighted pushforward density}\label{app:pushforward}

Adding \eqref{eq:char_ell} and \eqref{eq:app_liouville}, the divergence cancels:
\[
  \frac{\mathrm d}{\mathrm dt}\bigl(\ell_t+\log J_t\bigr)=g\bigl(t,X_t,e^{\ell_t}\bigr)=\dot\lw_t .
\]
Integrating from $0$ to $t$,
\begin{equation}\label{eq:app_jacobian}
  \mu\bigl(t,\Phi_t(x_0)\bigr)\,J_t(x_0)=\rho_0(x_0)\,e^{\lw_t(x_0)-\lw_0(x_0)} .
\end{equation}

\paragraph{Proof of Lemma~\ref{lem:weighted_residual}.}
(i) Equations \eqref{eq:char_ell} and \eqref{eq:char_score} and their initial values depend only on the starting point $X_0=x_0$, not on its distribution, so Sections~\ref{app:proof_carry} and~\ref{app:score} apply to every $x_0$.

(ii) With $\lw_0=\log\rho_0(x_0)-\log q_0(x_0)$, identity \eqref{eq:app_jacobian} reads $q_0(x_0)\,e^{\lw_t(x_0)}=\mu(t,\Phi_t(x_0))\,J_t(x_0)$. For such $f$, the change of variables $x=\Phi_t(x_0)$ gives
\begin{align*}
  \E\bigl[e^{\lw_t}f(t,X_t)\bigr]
  &=\int q_0(x_0)\,e^{\lw_t(x_0)}f\bigl(t,\Phi_t(x_0)\bigr)\dd x_0\\
  &=\int\mu\bigl(t,\Phi_t(x_0)\bigr)\,J_t(x_0)\,f\bigl(t,\Phi_t(x_0)\bigr)\dd x_0
  =\int\mu(t,x)\,f(t,x)\dd x .
\end{align*}
The integral over $x_0$ covers the support of $\rho_0$ because $q_0>0$ wherever $\rho_0>0$.

(iii) Take $f=1$ and $f=|R_v|^2$; by (i), $R_v(t,X_t)=v(t,X_t)-b(t,X_t)+\gamma\Sc_t$. \qed

\paragraph{Proof of Proposition~\ref{prop:carry}(ii).}
Take $q_0=\rho_0/M_0$, so that $\lw_0=\log M_0$, and $f=\mathbf 1_A$ in Lemma~\ref{lem:weighted_residual}(ii). This gives \eqref{eq:mu_def}; the existence and uniqueness of $\mu$ were shown in Section~\ref{app:proof_carry}. \qed

\subsection{Zero residual}\label{app:zero_residual}

\paragraph{Proof of Proposition~\ref{prop:selfconsistent}.}
If \eqref{eq:self_consistency} holds, then $\mu v=\mu b-\gamma\nabla\mu$, so
\[
  \nabla\!\cdot(\mu v)=\nabla\!\cdot(\mu b)-\gamma\Delta\mu ,
\]
and \eqref{eq:mu_pde} becomes \eqref{eq:rd_general} with $\mu(0)=\rho_0$. Since $\mu$ is a strictly positive classical solution, the uniqueness in Assumption~\ref{ass:continuum} gives $\mu=\rho$. Then $v=b-\gamma\nabla\log\rho=v_\rho$, and Proposition~\ref{prop:carry}(i) gives $\ell_t=\log\rho(t,X_t)$ and $\Sc_t=\nabla\log\rho(t,X_t)$.

Conversely, by \eqref{eq:rd_transport}, $\rho$ is a strictly positive classical solution of \eqref{eq:mu_pde} with $v=v_\rho$ and $\rho(0)=\rho_0$. If $v_\rho$ is admissible, the uniqueness in Section~\ref{app:proof_carry} gives $\mu[v_\rho]=\rho$, and hence $v_\rho=b-\gamma\nabla\log\mu[v_\rho]$. \qed

\subsection{The velocity equation}\label{app:proof_velocity}

\paragraph{Proof of Proposition~\ref{prop:velocity}.}
Let $v$ be $C^3$ in $x$ and $C^1$ in $t$ on each window, not necessarily admissible, let $u>0$ be a classical solution of $\partial_tu+\nabla\!\cdot(u\,v)=g(t,x,u)\,u$, and write $w=\nabla\log u$ and $R=v-b+\gamma w$. The computation of Section~\ref{app:score}, written in Eulerian form with $\nabla(v\cdot\nabla\log u)=\nabla(v\cdot w)$, gives
\[
  \partial_tw=-\nabla(v\cdot w)-\nabla(\nabla\!\cdot v)+\partial_xg+\partial_\ell g\,w ,
\]
with $g$ and its derivatives evaluated at $u$ and $\partial_\ell g=\partial_\rho g\,u$.
We replace $w$ by the residual. Since $\gamma w=R+b-v$, that is $w=R/\gamma+s_v$ with $s_v=(b-v)/\gamma$,
\[
  \gamma\nabla(v\cdot w)=\nabla(v\cdot R)+\nabla\bigl(v\cdot(b-v)\bigr),
  \qquad
  \gamma\bigl(\partial_xg+\partial_\ell g\,w\bigr)=\gamma\bigl(\partial_xg+\partial_\ell g\,s_v\bigr)+\partial_\ell g\,R .
\]
Substituting both into $\partial_tR=\partial_tv-\partial_tb+\gamma\,\partial_tw$ and collecting the terms that contain $R$ on the left,
\[
  \partial_tR+\nabla(v\cdot R)-\partial_\ell g\,R
  =\partial_tv-\partial_tb-\gamma\nabla(\nabla\!\cdot v)-\nabla\bigl(v\cdot(b-v)\bigr)+\gamma\bigl(\partial_xg+\partial_\ell g\,s_v\bigr).
\]
For $u=\mu[v]$, which solves \eqref{eq:mu_pde}, the right-hand side is $P(v)$ of \eqref{eq:velocity_pde_residual}, because $\gamma\,\partial_\ell g\,s_v=\mu\,\partial_\rho g\,(b-v)$, and the identity is \eqref{eq:residual_transport}. If $R_v\equiv0$ on $[0,T]$, the left-hand side vanishes, which is (i). For (ii), take $u=\rho$ and $v=v_\rho$: by \eqref{eq:rd_transport}, $\rho$ solves the same equation and $R=v_\rho-b+\gamma\nabla\log\rho=0$, so the right-hand side vanishes, which is \eqref{eq:velocity_equation}. This step uses only the equation for $\rho$ and not the flow of $v_\rho$, so it does not require $v_\rho$ to be admissible. \qed

\section{The relative-entropy estimate}\label{app:entropy}

This appendix proves Theorem~\ref{thm:stability}.

\subsection{Setting}\label{app:entropy_setting}

Let $v$ be admissible, $\mu=\mu[v]$ its induced density, and write
\[
  h:=\log\frac{\mu}{\rho},\qquad
  \Xi:=\rho\,\phi\Bigl(\frac{\mu}{\rho}\Bigr)=\mu h-\mu+\rho\;\ge0,\qquad
  \mathcal E(t):=\KL(\mu_t\|\rho_t)=\int\Xi\dd x,
\]
with $\phi(s)=s\log s-s+1$ as in \eqref{eq:gen_kl}, and
\[
  F(t):=\int\mu\,|R_v|^2\dd x ,\qquad g_\mu:=g(t,x,\mu(t,x)),\qquad g_\rho:=g(t,x,\rho(t,x)).
\]
We call $\Xi$ the KL density. By the definition of $R_v$, the transport-reaction equation \eqref{eq:mu_pde} and equation \eqref{eq:rd_general} read
\begin{align}
  \partial_t\mu+\nabla\!\cdot(\mu b)+\nabla\!\cdot(\mu R_v)&=\gamma\Delta\mu+g_\mu\mu ,\label{eq:st_perturbed}\\
  \partial_t\rho+\nabla\!\cdot(\rho b)&=\gamma\Delta\rho+g_\rho\rho .\label{eq:st_true}
\end{align}

The proof differentiates $\mathcal E$ under the integral and integrates by parts on $\R^d$, which requires the following assumption.

\begin{assumption}[Decay at infinity]\label{ass:decay}
For the admissible velocity $v$ under consideration, $\rho$ and $\mu[v]$ decay at infinity fast enough that $\mathcal E$ can be differentiated under the integral and the boundary terms produced by integration by parts on $\R^d$ vanish.
\end{assumption}

\subsection{The relative-entropy identity}\label{app:identity}

\begin{proposition}[Relative-entropy identity]\label{prop:st_identity}
Under Assumptions~\ref{ass:continuum} and~\ref{ass:decay}, for almost every $t$,
\begin{equation}\label{eq:st_identity}
  \mathcal E'(t)
  =\underbrace{-\gamma\int\mu|\nabla h|^2\dd x}_{\text{diffusion}}
   +\underbrace{\int\mu\,R_v\cdot\nabla h\dd x}_{\text{residual}}
   +\underbrace{\int g_\rho\,\Xi\dd x+\int(g_\mu-g_\rho)\,\mu h\dd x}_{\text{reaction}} .
\end{equation}
\end{proposition}

\begin{proof}
For $H(u,r)=u\log(u/r)-u+r$ one has $\partial_uH=\log(u/r)$ and $\partial_rH=1-u/r$, so
\[
  \mathcal E'(t)=\int\Bigl(h\,\partial_t\mu+\Bigl(1-\frac{\mu}{\rho}\Bigr)\partial_t\rho\Bigr)\dd x .
\]
We insert \eqref{eq:st_perturbed} and \eqref{eq:st_true} and treat the drift, diffusion, residual and reaction terms in turn, using
\begin{equation}\label{eq:app_two_gradients}
  \nabla\Bigl(1-\frac{\mu}{\rho}\Bigr)=-\frac{\mu}{\rho}\,\nabla h,
  \qquad
  \nabla\log\mu-\nabla\log\rho=\nabla h ,
\end{equation}
and integrating by parts under Assumption~\ref{ass:decay}.

\paragraph{Drift part.}
Integrating by parts and using \eqref{eq:app_two_gradients},
\[
  -\int h\,\nabla\!\cdot(\mu b)\dd x-\int\Bigl(1-\frac{\mu}{\rho}\Bigr)\nabla\!\cdot(\rho b)\dd x
  =\int\mu\,b\cdot\nabla h\dd x-\int\rho\,b\cdot\frac{\mu}{\rho}\nabla h\dd x=0 .
\]

\paragraph{Diffusion part.}
Integrating by parts and using \eqref{eq:app_two_gradients},
\begin{align*}
  \gamma\int h\,\Delta\mu\dd x+\gamma\int\Bigl(1-\frac{\mu}{\rho}\Bigr)\Delta\rho\dd x
  &=-\gamma\int\nabla h\cdot\nabla\mu\dd x+\gamma\int\frac{\mu}{\rho}\,\nabla h\cdot\nabla\rho\dd x\\
  &=-\gamma\int\mu\,\nabla h\cdot\bigl(\nabla\log\mu-\nabla\log\rho\bigr)\dd x
  =-\gamma\int\mu|\nabla h|^2\dd x .
\end{align*}

\paragraph{Residual part.}
Only \eqref{eq:st_perturbed} contains $R_v$, and
\[
  -\int h\,\nabla\!\cdot(\mu R_v)\dd x=\int\mu\,R_v\cdot\nabla h\dd x .
\]

\paragraph{Reaction part.}
Adding and subtracting $g_\rho\mu h$,
\[
  \int g_\mu\,\mu h\dd x+\int g_\rho\,(\rho-\mu)\dd x
  =\int g_\rho\,(\mu h-\mu+\rho)\dd x+\int(g_\mu-g_\rho)\,\mu h\dd x
  =\int g_\rho\,\Xi\dd x+\int(g_\mu-g_\rho)\,\mu h\dd x .
\]

Adding the four parts gives \eqref{eq:st_identity}.
\end{proof}

\subsection{From the identity to the residual inequality}\label{app:inequality}

Since $g$ is non-increasing in the density, $g_\mu-g_\rho\le0$ where $\mu>\rho$, that is where $h>0$, and $g_\mu-g_\rho\ge0$ where $h<0$. Hence $(g_\mu-g_\rho)\,h\le0$ pointwise, and the last term of \eqref{eq:st_identity} is non-positive. For the WFR gradient flow, for example, $g_\mu-g_\rho=-\beta h$, so this term is $-\beta\int\mu h^2\dd x$.

The diffusion and residual terms combine pointwise as
\[
  -\gamma|\nabla h|^2+R_v\cdot\nabla h
  =-\gamma\Bigl|\nabla h-\frac{R_v}{2\gamma}\Bigr|^2+\frac{|R_v|^2}{4\gamma}
  \le\frac{|R_v|^2}{4\gamma}.
\]
Multiplying by $\mu$, integrating and dropping the non-positive reaction term gives
\begin{equation}\label{eq:st_basic}
  \mathcal E'(t)\le\frac{F(t)}{4\gamma}+\int g_\rho\,\Xi\dd x\le\frac{F(t)}{4\gamma}+\int[g_\rho]_+\,\Xi\dd x .
\end{equation}
The first inequality is Theorem~\ref{thm:stability}(i).

\subsection{Proof of Theorem~\ref{thm:stability}(ii) and the Fisher--KPP constant}\label{app:stability}

\begin{proof}[Proof of Theorem~\ref{thm:stability}(ii)]
By the definition of $G$, $[g_\rho]_+\le G(t)$, so \eqref{eq:st_basic} gives
\[
  \mathcal E'(t)\le G(t)\,\mathcal E(t)+\frac{F(t)}{4\gamma}.
\]
Multiplying by $e^{-\int_0^tG(\tau)\,\mathrm d\tau}$ gives $\frac{\mathrm d}{\mathrm dt}\bigl(e^{-\int_0^tG(\tau)\,\mathrm d\tau}\mathcal E(t)\bigr)\le e^{-\int_0^tG(\tau)\,\mathrm d\tau}F(t)/(4\gamma)$. Integrating from $0$ to $t$ with $\mathcal E(0)=0$,
\[
  \mathcal E(t)\le\frac1{4\gamma}\int_0^te^{\int_s^tG(\tau)\,\mathrm d\tau}F(s)\dd s
  =\frac1{4\gamma}\int_0^tM_\mu(s)\,e^{\int_s^tG(\tau)\,\mathrm d\tau}\,\frac{F(s)}{M_\mu(s)}\dd s
  \le\frac{C_G}{4\gamma}\,J_{\mathrm{reg}}(v),
\]
where the last step bounds $M_\mu(s)e^{\int_s^tG(\tau)\,\mathrm d\tau}$ by $C_G$ and the remaining integral by $J_{\mathrm{reg}}(v)$.
\end{proof}

For the Fisher--KPP reaction of Section~\ref{sec:kpp}, $g(t,x,u)=k(1-u/K)\le k$ for every $u\ge0$. Hence $G\le k$, and also $M_\mu'=\int g(t,x,\mu)\,\mu\dd x\le kM_\mu$, so $M_\mu(s)\le M_0e^{ks}$ and $C_G\le M_0e^{kT}$. Since the KL divergence scales with the mass, the bound reads, relative to the initial mass,
\[
  \sup_{0\le t\le T}\frac{\KL(\mu_t\|\rho_t)}{M_0}\le\frac{e^{kT}}{4\gamma}\,J_{\mathrm{reg}}(v),
\]
with $e^{kT}=e^2\approx7.4$ for the one-dimensional problem ($k=1$, $T=2$) and $e^4\approx55$ for the three-dimensional one ($k=2$, $T=2$).

When $[g_\rho]_+$ is unbounded in $x$, as for the linear source of Section~\ref{sec:ou} and for the WFR gradient flow, $G\equiv\infty$ and part (ii) does not apply. Closing \eqref{eq:st_basic} then requires control of the tails of $\Xi$ that holds uniformly over the candidate velocities; we have not established such control for the velocities produced by Algorithm~\ref{alg:main}.

\section{Implementation}\label{app:impl}

\subsection{Velocity representation}\label{sec:coords}
All initial densities in Section~\ref{sec:numerics} have the form $\rho_0=M_0\sum_k w_k\,\mathcal N(m_k,\Sigma_k)$, a Gaussian mixture multiplied by the mass $M_0$. Their heat flow
\begin{equation}\label{eq:heat_reference}
  \refeta(t,x)=M_0\sum_k w_k\,\mathcal N\bigl(x;m_k,\Sigma_k+2\gamma tI\bigr)
\end{equation}
solves $\partial_t\refeta=\gamma\Delta\refeta$ with $\refeta(0)=\rho_0$, and it and its derivatives are available in closed form. We write the velocity as
\begin{equation}\label{eq:velocity_split}
  v_\theta=\veta+\uth,\qquad \veta:=-\gamma\nabla\log\refeta ,
\end{equation}
so that the network represents only the deviation from the pure heat flow, which is caused by the drift and the reaction. The network $\uth$ is a multilayer perceptron with SiLU activations~\cite{elfwing2018silu} in double precision, and its last layer is initialized to zero. Its output is bounded componentwise through a $\tanh$, and its spatial input is scaled by the width of the heat flow; Table~\ref{tab:settings} gives the sizes.

Instead of $\ell$ we carry $\cc:=\ell-\log\refeta(t,X)$, which vanishes at $t=0$ and stays zero along the pure heat flow. By the log-density equation \eqref{eq:char_ell} and the heat equation for $\refeta$, it satisfies $\dot\cc=g(t,X,\refeta e^{\cc})-\nabla\cdot\uth-\uth\cdot\nabla\log\refeta$, and Heun's method \eqref{eq:heun} is applied to $(X,\cc,\Sc,\lw)$. The score $\Sc$ is carried in full by the score equation \eqref{eq:char_score} with $v=\veta+\uth$; the derivatives of $\veta$ are analytic, and those of $\uth$ ($(\nabla\uth)^\top\Sc$, $\nabla\cdot\uth$ and $\nabla(\nabla\cdot\uth)$) are computed by automatic differentiation. The regression labels are $\lab=b-\gamma\Sc-\veta$, so that the residual \eqref{eq:residual} is $R_v=\uth-\lab$.

\subsection{Training loss and stopping rule}\label{app:iteration}

Both terms of \eqref{eq:J_modular} use the Huber function with threshold one (the smooth $L^1$ loss), applied to each component and summed,
\[
  \psi(r)=\sum_{c=1}^dH(r_c),\qquad
  H(s)=\begin{cases}s^2/2,&|s|<1,\\|s|-1/2,&\text{otherwise}.\end{cases}
\]
Its derivative is bounded by one in each component, so the derivative of the loss with respect to each residual component is bounded by the weight of the particle. The regression term uses $\psi_{\rm reg}(r)=\psi(r)$. The regularizer term applies $\psi$ to the residual \eqref{eq:P_theta} divided componentwise by one plus the sum of the magnitudes of its terms.

In each window the loss \eqref{eq:J_modular} is minimized with the adaptive weights, and each optimizer step uses particles at one randomly chosen time $t_n$. Before window 0 the network is pretrained on the exact velocity $b-\gamma\nabla\log\rho_0$ at $t=0$, and during window 0 this pretraining loss stays in the objective with a small weight. An outer iteration is accepted when the relative change of the labels and the relative fit residual are below fixed thresholds. A window stops after two consecutive accepted iterations or at the maximum number of iterations. Table~\ref{tab:settings} lists the batch sizes, the thresholds and the iteration limits.

\subsection{Period-adapted sampling and adaptive weights}\label{app:banks}

The period-adapted sampling distributions and the adaptive weights both need a rough prediction of where the mass will be. They use a predicted density $\rho_{\mathrm{pred}}$ built from the problem data only, which is different from the heat flow $\refeta$ of Appendix~\ref{sec:coords} and from the reference solutions of Appendix~\ref{app:evaluation}; Table~\ref{tab:banks} lists it, together with the defensive density used below, for each problem. Each $\rho_{\mathrm{pred}}$ is built in closed form from simplified problems: OU without the source, diffusion and reaction applied one after the other for Fisher--KPP, and one target component at a time for WFR. It enters only the sampling distributions and the adaptive weights; the labels and the residual \eqref{eq:residual} do not use it.

\begin{table}[ht]
\centering
\caption{Problem-dependent quantities of the period-adapted sampling.}
\label{tab:banks}
\smallskip
\small
\begin{tabular}{@{}lp{0.45\textwidth}p{0.3\textwidth}@{}}
\toprule
Problem & Predicted density $\rho_{\mathrm{pred}}$ & Defensive density $\pi_{\mathrm d}$ \\
\midrule
OU & solution without source & invariant density of the drift \\
Fisher--KPP & heat flow followed by the logistic map, $\refeta/\bigl(e^{-kt}+(1-e^{-kt})\refeta/K\bigr)$ & $\mathcal N(0,9)$ \\
WFR & mixture of the exact WFR gradient flows from the Gaussian $\rho_0$ towards each single target component & the target $\pi$ \\
\bottomrule
\end{tabular}
\end{table}

 The particles are moved by $v_{\mathrm{pred}}=b-\gamma\nabla\log\rho_{\mathrm{pred}}$ with Heun's method and step $\Delta t$, and their log-weights $\widehat\lw_i(t)$ accumulate $g$ evaluated at $\rho_{\mathrm{pred}}$ (for OU, the actual source $a\cdot x$). The predicted growth $\Ghat_i^{(j)}$ in \eqref{eq:measure} is computed before training and kept fixed.

For each period $\mathcal T_m$, the sampling distribution $q_0^m$ should place particles where the mass of the period comes from. We integrate pilot particles under the predicted flow and weight each initial position by the particle's share of the mass at several times in the period; for a single time, this weighted distribution is the importance sampling density with zero variance for the mass carried by the predicted flow. A Gaussian mixture $q_{m,k}$ is fitted to these weighted positions separately for each mode $k$ (a target component for WFR; OU and Fisher--KPP have a single mode), and
\[
  q_0^m=(1-\epsilon_{\mathrm d})\sum_k p_k\,q_{m,k}+\epsilon_{\mathrm d}\,q_0^{\mathrm c},\qquad p_k\propto\max(f_{\min},f_k),\qquad q_0^{\mathrm c}=(1-\epsilon_{\mathrm c})\,\rho_0/M_0+\epsilon_{\mathrm c}\,\pi_{\mathrm d},
\]
where $f_k$ is the share of mode $k$ in the period and $\pi_{\mathrm d}$ is given in Table~\ref{tab:banks}. The defensive mixture $q_0^{\mathrm c}$ keeps $q_0^m>0$ wherever $\rho_0>0$, as Lemma~\ref{lem:weighted_residual} requires, and the lower bound $f_{\min}$ keeps particles in every mode. The fit is repeated once with particles drawn from the first $q_0^m$. The single distribution in the comparison of Section~\ref{sec:wide} is $q_0^{\mathrm c}$. Table~\ref{tab:settings} lists the constants.

\subsection{Solver settings}\label{app:settings}
Table~\ref{tab:settings} lists the settings and all numerical constants of Appendices~\ref{sec:coords}--\ref{app:banks}. The problem-dependent quantities are the number of windows $T/W$, Table~\ref{tab:banks} and the evaluation set-up of Appendix~\ref{app:evaluation}.

\begin{table}[ht]
\centering
\caption{Solver settings. Entries are shared by all problems unless marked.}
\label{tab:settings}
\smallskip
\small
\begin{tabular}{@{}lp{0.72\textwidth}@{}}
\toprule
Component & Setting \\
\midrule
Velocity & $v_\theta=\veta+\uth$ (Appendix~\ref{sec:coords}); $\uth$: 3 hidden layers of width 256, SiLU, double precision, output $R_u\tanh(\cdot/R_u)$ with $R_u=640$ \\
Integration & Heun, $\Delta t=0.0025$; windows of width $W=0.01$ (100 windows for $T=1$, 200 for $T=2$) \\
Periods & $\mathcal T_0=[0,0.125T)$, $\mathcal T_1=[0.125T,0.4T)$, $\mathcal T_2=[0.4T,T]$ \\
Sampling & three period-adapted sampling distributions with 8192 particles each (32768 for 3D Fisher--KPP); pilot of 262144 initial points; second round of 131072 per distribution; at most 4 Gaussian components per mode (8 for 3D Fisher--KPP); $\epsilon_{\mathrm c}=0.375$, $\epsilon_{\mathrm d}=0.2$, $f_{\min}=0.1$ \\
Loss & Huber (threshold 1), $\lamres=0.1$, adaptive weights \eqref{eq:measure} \\
Optimizer & Adam~\cite{kingma2015adam}, learning rate $3\times10^{-4}$, minibatch 2048 (1024 for the regularizer term), gradient-norm clip 10; Adam is reset and the network warm-started at each window \\
Outer iterations & 150 updates per iteration; an iteration is accepted when the relative label change is at most 0.01 and the relative fit residual at most 0.05; a window stops after two consecutive accepted iterations, with 4 to 40 iterations per window \\
Window 0 & 1000 pretraining updates; weight 0.2 of the pretraining loss \\
Seeds & three training seeds; one A100 (40\,GB) per run \\
\bottomrule
\end{tabular}
\end{table}

\section{Reference solutions and evaluation}\label{app:evaluation}

\subsection{Exact solution of the OU equation with a linear source}\label{app:ou_exact}

\begin{proposition}\label{prop:ou_exact}
Let $d\ge1$, $\kappa,\gamma>0$, $a,m_0\in\mathbb R^d$, $M_0>0$, and let $C_0$ be symmetric positive definite. Then $\rho(t)=M(t)\,\mathcal N\bigl(m(t),C(t)\bigr)$, with $C$, $m$ and $M$ given by \eqref{eq:ou_exact}, solves the OU equation \eqref{eq:ou_pde},
\[
  \partial_t\rho=\nabla\cdot(\kappa x\,\rho)+\gamma\Delta\rho+(a^\top x)\,\rho,\qquad \rho(0)=M_0\,\mathcal N(m_0,C_0),
\]
and its score and velocity \eqref{eq:ou_fields} are
\[
  \nabla\log\rho(t,x)=-C(t)^{-1}\bigl(x-m(t)\bigr),\qquad v_\rho(t,x)=-\kappa x+\gamma C(t)^{-1}\bigl(x-m(t)\bigr).
\]
\end{proposition}

\begin{proof}
Consider the Gaussian form \eqref{eq:ou_density} with unknown functions $M(t)>0$, $m(t)$ and $C(t)>0$. Set $y=x-m(t)$ and $Q=C(t)^{-1}$. Differentiation gives
\begin{align}
  \frac{\partial_t\rho}{\rho}
    &=\frac{\dot M}{M}-\frac12\operatorname{tr}(Q\dot C)
      +\dot m^\top Qy+\frac12y^\top Q\dot C Qy,
      \label{eq:ou_time_derivative}\\
  \frac{\nabla\rho}{\rho}&=-Qy,
  \qquad
  \frac{\Delta\rho}{\rho}=y^\top Q^2y-\operatorname{tr}Q.
      \label{eq:ou_spatial_derivatives}
\end{align}
The right-hand side of \eqref{eq:ou_pde}, divided by $\rho$, is therefore
\begin{equation}\label{eq:ou_pde_expansion}
  \begin{aligned}
    \frac{\nabla\!\cdot(\kappa x\rho)+\gamma\Delta\rho
          +(a^\top x)\rho}{\rho}
      ={}&y^\top(-\kappa Q+\gamma Q^2)y
         +(-\kappa m+Ca)^\top Qy\\
        &+\kappa d-\gamma\operatorname{tr}Q+a^\top m.
  \end{aligned}
\end{equation}
Equating the quadratic, linear and constant coefficients in $y$ with those in \eqref{eq:ou_time_derivative} gives
\begin{equation}\label{eq:ou_parameter_odes}
  \dot C=-2\kappa C+2\gamma I,
  \qquad
  \dot m=-\kappa m+Ca,
  \qquad
  \frac{\mathrm d}{\mathrm dt}\log M=a^\top m,
\end{equation}
with $C(0)=C_0$, $m(0)=m_0$ and $M(0)=M_0$. For the constant coefficient we used $-\tfrac12\operatorname{tr}(Q\dot C)=\kappa d-\gamma\operatorname{tr}Q$.

Solving the covariance equation gives
\[
  C(t)=e^{-2\kappa t}C_0+\frac{\gamma}{\kappa}(1-e^{-2\kappa t})I.
\]
The integrating factor for the mean equation then yields
\[
  m(t)=e^{-\kappa t}m_0+\int_0^t e^{-\kappa(t-s)}C(s)a\,\mathrm ds
       =r_t m_0+\frac{r_t(1-r_t)}{\kappa}C_0a+\frac{\gamma}{\kappa^2}(1-r_t)^2a.
\]
Finally,
\begin{equation*}
  \begin{aligned}
    \log\frac{M(t)}{M_0}
      &=\int_0^t a^\top m(s)\,\mathrm ds\\
      &=\frac{1-r_t}{\kappa}a^\top m_0
        +\frac{(1-r_t)^2}{2\kappa^2}a^\top C_0a
        +\frac{\gamma}{2\kappa^3}(2\kappa t-3+4r_t-r_t^2)|a|^2,
  \end{aligned}
\end{equation*}
where we used
\[
  \int_0^t(e^{-\kappa s}-e^{-2\kappa s})\,\mathrm ds=\frac{(1-r_t)^2}{2\kappa},
  \qquad
  \int_0^t(1-e^{-\kappa s})^2\,\mathrm ds=\frac{2\kappa t-3+4r_t-r_t^2}{2\kappa}.
\]
These are the formulas in \eqref{eq:ou_exact}. The matrix $C(t)$ remains positive definite and $M(t)$ is positive and finite at every finite time, so \eqref{eq:ou_density} defines a smooth, positive density of finite mass. It satisfies \eqref{eq:ou_pde} and the prescribed initial condition; uniqueness under Assumption~\ref{ass:continuum} identifies it with the reference solution. Equation~\eqref{eq:ou_fields} follows from $\nabla\log\rho=-Qy$ and $v_\rho=b-\gamma\nabla\log\rho$.
\end{proof}

\subsection{Numerical reference solutions}\label{app:references}
The reference solutions use finite volumes. Fisher--KPP in one dimension uses $[-25,25]$ with 1025 and 2049 cells and Strang splitting~\cite{strang1968construction}, with an exact logistic half-step on each side of a Crank--Nicolson diffusion step~\cite{crank1947practical}. Fisher--KPP in three dimensions solves the radial equation on $r\in[0,25]$ with 513 and 1025 points. WFR uses $[-17.5,17.5]^2$ with $2049^2$ and $4097^2$ points, exact reaction half-steps and Chang--Cooper fluxes~\cite{chang1970practical}. In the relative density and score, the two resolutions differ by at most $2.2\times10^{-4}$ for Fisher--KPP in one dimension and $4.6\times10^{-4}$ in three dimensions, and by at most $3.0\times10^{-5}$ in the mode fractions of WFR.

\subsection{Output particles and metrics}\label{app:metrics}
After training, the window fields are frozen and integrated from new initial points drawn from the equal mixture of the three period-adapted sampling distributions, with two independent output seeds: $N=32768$ particles for $d\le2$ and $65536$ for $d=3$. These particles are never used for training, and in the comparison of Section~\ref{sec:wide} all runs are evaluated on the same output particles. The report times are given in Section~\ref{sec:numerics}; for Fisher--KPP they are $\{0.02,0.05,0.1,0.2,0.4,0.8,1.2,1.6,2\}$.

Let $\bar q_0$ be this mixture. The particles start from $\lambda_i(0)=\log\rho_0(X_i(0))-\log\bar q_0(X_i(0))$, so the sums below estimate expectations over $X_0\sim\bar q_0$. At a report time $t$, let $w_i=e^{\lambda_i(t)}$ and $\bar w_i=w_i/\sum_r w_r$.

\paragraph{Mass.}
\[
  \widehat M(t)=\frac1N\sum_iw_i,\qquad
  e_M(t)=\left|\frac{\widehat M(t)}{M(t)}-1\right|,
  \qquad M(t)=\int\rho(t,x)\,\mathrm dx.
\]

\paragraph{Distribution.}
For a partition of unity $\sum_c\chi_c(x)=1$,
\[
  \widehat f_c(t)=\sum_i\bar w_i\,\chi_c(X_i(t)),\qquad
  f_c(t)=\frac1{M(t)}\int\rho(t,x)\chi_c(x)\,\mathrm dx,\qquad
  D(t)=\sum_c|\widehat f_c(t)-f_c(t)|.
\]
For OU and Fisher--KPP, $\chi_c$ are cell indicators: $16\times16$ cells on $[-25,25]^2$ for OU, 128 cells on $[-25,25]$ and $12^3$ cells on $[-18,18]^3$ for Fisher--KPP, each with an additional cell for the outside of the box. For the WFR problem they are the target-component responsibilities $r_k(x)=\tfrac14\mathcal N(x;m_k,\Sigma_k)/\pi(x)$.

\paragraph{Velocity and score.}
On independent points drawn from $\rho(t)/M(t)$, let $e_v(t)^2$ and $a_v(t)^2$ be the sample means of $|v_\theta-v_\rho|^2$ and $|v_\rho|^2$. With $\mathcal T$ the positive report times,
\[
  E_v=\left(\frac{\sum_{t\in\mathcal T}e_v(t)^2}{\sum_{t\in\mathcal T}a_v(t)^2}\right)^{1/2},
  \qquad
  E_S(t)=\left(\frac{\sum_i\bar w_i|S_i-\nabla\log\rho(t,X_i)|^2}
                         {\sum_i\bar w_i|\nabla\log\rho(t,X_i)|^2}\right)^{1/2}.
\]
 In the WFR problem, particles outside the reference grid are excluded from the pointwise comparisons. For OU, $\nabla v_\rho(t)=-\kappa I+\gamma C(t)^{-1}$ is constant in space and $\nabla(\nabla\cdot v_\rho)=0$. Section~\ref{sec:ou} reports the absolute root-mean-square errors
\[
  \Bigl(\frac1n\sum_{k=1}^n\bigl|\nabla v_\theta(t,y_k)-\nabla v_\rho(t)\bigr|_F^2\Bigr)^{1/2},\qquad
  \Bigl(\frac1n\sum_{k=1}^n\bigl|\nabla(\nabla\cdot v_\theta)(t,y_k)\bigr|^2\Bigr)^{1/2}
\]
on $n=128$ independent points $y_k$ drawn from $\rho(t)/M(t)$, maximized over the training seeds.

\paragraph{Figures.} Snapshots and profiles use one training seed and one output seed. Particle histograms divide the mass $\frac1N\sum_i w_i$ of the particles in each bin of the plotted coordinate by the bin width, without smoothing or renormalization.

\end{document}